\documentclass{amsart}

\usepackage{tikz}

\title{Maximizers for the {S}trichartz inequality}

\author[D.\ Foschi]%
{Damiano Foschi}

\date{June 9, 2006}


\newcommand{\abs}[1]{\left|#1\right|}
\newcommand{\Abs}[1]{\bigl|#1\bigr|}
\newcommand{\norm}[1]{\left\|#1\right\|}
\newcommand{\Norm}[1]{\bigl\|#1\bigr\|}
\newcommand{\tonde}[1]{\left(#1\right)}
\newcommand{\Tonde}[1]{\bigl(#1\bigr)}
\newcommand{\TOnde}[1]{\Bigl(#1\Bigr)}

\newcommand{\graffe}[1]{\left\{#1\right\}}
\newcommand{\Graffe}[1]{\bigl\{#1\bigr\}}
\newcommand{\ip}[1]{\langle #1 \rangle}
\newcommand{\pvector}[1]{\begin{pmatrix}#1\end{pmatrix}}
\newcommand{\ddirac}[1]{%
  \,\boldsymbol{\delta}\!\pvector{#1}\!}
\newcommand{\MU}{\boldsymbol{\mu}}
\newcommand{\chic}{\boldsymbol{\chi}}
\newcommand{\R}{\mathbb{R}}  
\newcommand{\C}{\mathbb{C}}  
\newcommand{\iC}{\mathrm{i}} 
\newcommand{\Eu}{\mathrm{e}} 
\renewcommand{\Re}{{\mathfrak{Re}}}
\renewcommand{\Im}{{\mathfrak{Im}}}
\newcommand{\intersection}{\cap}
\newcommand{\union}{\cup} 

\newcommand{\mC}{{\mathcal C}}
\newcommand{\mE}{{\mathcal E}}
\newcommand{\mG}{{\mathcal G}}
\newcommand{\mH}{{\mathcal H}}
\newcommand{\mL}{{\mathcal L}}
\newcommand{\mN}{{\mathcal N}}
\newcommand{\mP}{{\mathcal P}}
\newcommand{\de}{\partial}
\newcommand{\pder}[3][]{{\frac{\de^{#1} #2}{\de{#3}^{#1}}}}
\renewcommand{\d}{\,{\rm d}}
\newcommand{\1}[1]{\frac{1}{#1}}
\newcommand{\conv}{\ast}
\newcommand{\ovl}{\overline}
\newcommand{\changeto}{\rightsquigarrow}
\newcommand{\ut}{\widetilde{u}}
\newcommand{\Ft}{\widetilde{F}}
\newcommand{\Gt}{\widetilde{G}}
\newcommand{\fh}{\widehat{f}}
\newcommand{\fph}{\widehat{f_+}}
\newcommand{\fmh}{\widehat{f_-}}
\newcommand{\Ah}{\widehat{A}}
\newcommand{\bh}{\widehat{b}}
\newcommand{\Bh}{\widehat{B}}
\newcommand{\Ch}{\widehat{C}}

\theoremstyle{plain}
\newtheorem{theorem}{Theorem}[section]
\newtheorem{lemma}[theorem]{Lemma}
\newtheorem{proposition}[theorem]{Proposition}
\newtheorem{conjecture}[theorem]{Conjecture}
\theoremstyle{definition}
\newtheorem{definition}[theorem]{Definition}
\theoremstyle{remark}
\newtheorem{remark}[theorem]{Remark}

\begingroup\makeatletter\ifx\SetFigFont\undefined%
\gdef\SetFigFont#1#2#3#4#5{%
  \reset@font\fontsize{#1}{#2pt}%
  \fontfamily{#3}\fontseries{#4}\fontshape{#5}%
  \selectfont}%
\fi\endgroup%

\begin{document}

\maketitle

\begin{center}
  Dipartimento di Matematica \\
  Universit\`a di Ferrara \\
  via Macchiavelli 34, I-44100 Ferrara - ITALY \\
  email: \verb|damiano.foschi@unife.it|
\end{center}

\begin{abstract}
  We compute explicitly the best constants
  and, by solving some functional equations,
  we find all maximizers for homogeneous Strichartz estimates
  for the Schr\"odinger equation and for the wave equation
  in the cases when the Lebesgue exponent is an even integer.
\end{abstract}

\section{%
  Introduction}
\label{sec:introduction}

Let $n$ be a positive integer and let $p = p(n) = 2 + 4/n$.
The Strichartz inequality
for the homogeneous Schr\"odinger equation
in $n$ spatial dimensions
states that there exists a constant $S > 0$ such that
\begin{equation} \label{eq:1}
  \norm{u}_{L^{p(n)}(\R^{1 + n})} \le S \norm{f}_{L^2(\R^n)},
\end{equation}
whenever $u(t, x)$ is the solution of the equation
\begin{equation} \label{eq:2}
  \iC \de_t u = \Delta u,
\end{equation}
with initial data $u(0, x) = f(x)$;
see~\cite{Str1977} for the original proof by Strichartz.
We denote by $S(n)$ the best constant for the estimate~\eqref{eq:1},
\begin{equation*}
  S(n) = \sup_{f \in L^2(\R^n)}
  \frac{\norm{\Eu^{-\iC t \Delta} f}_{L^{p(n)}(\R^{1 + n})}}
  {\norm{f}_{L^2(\R^n)}}.
\end{equation*}

If $n \ge 2$, we can also consider the Strichartz inequality
for the homogeneous wave equation
in $n$ spatial dimensions which
states that there exists a constant $W > 0$ such that
\begin{equation} \label{eq:3}
  \norm{u}_{L^{p(n - 1)}(\R^{1 + n})} \le
  W \Norm{(f, g)}_{\dot{H}^{\12}(\R^n)
    \times \dot{H}^{-\12}(\R^n)},
\end{equation}
whenever $u(t, x)$ is the solution of the equation
\begin{equation} \label{eq:4}
  \de_t^2 u = \Delta u,
\end{equation}
with initial data
\begin{equation} \label{eq:5}
  u(0, x) = f(x), \qquad \de_t u(0, x) = g(x).
\end{equation}
This also was proved in~\cite{Str1977}.
We denote by $W(n)$
the best constant for the estimate~\eqref{eq:3},
\begin{equation*}
  W(n) = \sup_{\substack{f \in \dot{H}^{\12}(\R^n) \\
      g \in \dot{H}^{-\12}(\R^n)}}
  \frac{\norm{\cos(t \sqrt{-\Delta}) f +
      \frac{\sin(t \sqrt{-\Delta})}{\sqrt{-\Delta}} g}_{%
      L^{p(n - 1)}(\R^{1 + n})}}
  {\norm{(f, g)}_{\dot{H}^{\12}(\R^n)
      \times \dot{H}^{-\12}(\R^n)}}.
\end{equation*}

Kunze~\cite{Kun2003} has recently proved
the existence of a maximizing function $f_* \in L^2(\R)$
for the estimate~\eqref{eq:1} in the special case $n = 1$ and $p = 6$,
which means that for the corresponding solution $u_*$,
\begin{equation*}
  \iC \de_t u_* = \de_x^2 u_*, \qquad u_*(0, x) = f_*(x),
\end{equation*}
we have the equality
$\norm{u_*}_{L^6(\R \times \R)} = S(1) \norm{f_*}_{L^2(\R)}$.
The proof in~\cite{Kun2003} is based on an elaborate application
of the concentration compactness principle
and does not provide an explicit expression for a maximizer.

Here, we present a more direct and elementary approach
which allows us to explicitly determine the families of maximizers
and compute the best constants
for the estimates~\eqref{eq:1} and~\eqref{eq:3}
when the exponent $p = p(n)$ is an even integer.
We show that the classes of maximizers are unique modulo
the natural geometric invariance properties of the equations.
Moreover, maximizers turn out to be smooth solutions
to some functional equations which can be solved explicitly.

\medskip
For the Schr\"odinger equation we have:

\begin{theorem}
  In the case $n = 1$ and $p = 6$, we have $S(1) = {12}^{-1/12}$;
  in the case $n = 2$ and $p = 4$, we have $S(2) = 2^{-1/2}$.
  In both cases an example of a maximizer $f_* \in L^2(\R^n)$
  for which we have
  \begin{equation} \label{eq:6}
    \norm{\Eu^{-\iC t \Delta} f_*}_{L^p(\R \times \R^n)} =
    S(n) \norm{f_*}_{L^2(\R^n)},
  \end{equation}
  is provided by the Gaussian function
  $f_*(x) = \exp\Tonde{-\abs{x}^2}$.
\end{theorem}

The geometric invariance properties of the equation~\eqref{eq:2}
suggest a way to completely characterize the class of all maximizers.

\begin{definition} \label{def:1}
  Let $\mG$ be the Lie group of transformations generated by:
  \begin{itemize}
  \item space-time translations:
    $u(t, x) \changeto u(t + t_0, x + x_0)$,
    with $t_0 \in \R$, $x_0 \in \R^n$;
  \item parabolic dilations:
    $u(t, x) \changeto u(\lambda^2 t, \lambda x)$,
    with $\lambda > 0$;
  \item change of scale:
    $u(t, x) \changeto \mu u(t, x)$, with $\mu > 0$;
  \item space rotations:
      $u(t, x) \changeto u(t, R x)$, 
      with $R \in SO(n)$;
  \item phase shifts:
      $u(t, x) \changeto \Eu^{\iC \theta} u(t, x)$, 
      with $\theta \in \R$;
  \item Galilean transformations:
    \begin{equation*}
      u(t, x) \changeto
      \exp\TOnde{\frac{\iC}4 \Tonde{\abs{v}^2 t + 2 v \cdot x}}
      u(t, x + t v),
    \end{equation*}
    with $v \in \R^n$.
  \end{itemize}
\end{definition}

If $u$ solves~\eqref{eq:2} and $g \in \mG$
then $v = g \cdot u$ is still a solution to~\eqref{eq:2}.
Moreover, the ratio $\norm{u}_{L^{p(n)}} / \norm{u(0)}_{L^2}$
is left unchanged by the action of $\mG$.

\begin{remark}
  We should mention that
  there exists another important (discrete) symmetry
  for the Schr\"odinger equation
  given by the \emph{pseudo-conformal inversion}:
  \begin{equation*}
    u(t, x) \changeto 
    t^{-n / 2} \Eu^{\iC \abs{x}^2 / (4 t)} u\TOnde{-\1t, \frac{x}{t}}.
  \end{equation*}
  Combining the inversion with translations and dilations,
  we obtain that the Schr\"odinger equation is invariant
  under the representation of $SL(2,\R)$ given by
  \begin{equation*}
    u(t, x) \changeto 
    (a + b t)^{-n / 2} \Eu^{\iC b \abs{x}^2 / (4 (a + b t))}
    u\TOnde{\frac{c + d t}{a + b t}, \frac{x}{a + b t}}, \qquad
    a d - b c = 1.
  \end{equation*}
  These transformations have many important applications.
  However, we do not really need them in the context of our analysis
  and for simplicity we are not including them
  in the list of generators of the group $\mG$.
\end{remark}

\begin{theorem}
  Let $(n, p) = (1, 6)$ or $(n, p) = (2, 4)$.
  Let $f_*(x) = \exp\Tonde{-\abs{x}^2}$ 
  and $u_*(t, x) = \Eu^{-\iC t \Delta} f_*(x)$
  be the corresponding solution
  to the Schr\"odinger equation~\eqref{eq:2}.
  Then the set of maximizers
  for which the equality~\eqref{eq:6} holds
  coincides with the set of initial data
  of solutions to~\eqref{eq:2} in the orbit of $u_*$
  under the action of the group $\mG$.
  In particular, all maximizers
  are given by $L^2$ functions of the form
  \begin{equation*}
    f_*(x) = \exp\Tonde{A \abs{x}^2 + b \cdot x + C},
  \end{equation*}
  with $A, C \in \C$, $b \in \C^n$ and $\Re(A) < 0$.
\end{theorem}

\medskip
For the wave equation we have:

\begin{theorem}
  In the case $n = 2$ and $p = 6$,
  we have $W(2) = (25 / 64 \pi)^{1/6}$;
  in the case $n = 3$ and $p = 4$,
  we have $W(3) = (3 / 16\pi)^{1/4}$.
  In both cases an example of a maximizer pair
  $(f_*, g_*) \in \dot{H}^{\12}(\R^n) \times \dot{H}^{-\12}(\R^n)$
  for which we have
  \begin{equation} \label{eq:7}
    \norm{\cos(t \sqrt{-\Delta}) f_* +
      \frac{\sin(t \sqrt{-\Delta})}{\sqrt{-\Delta}} g_*}_{%
      L^p(\R^{1 + n})} =
    W(n) \Norm{(f_*, g_*)}_{\dot{H}^{\12}(\R^n)
      \times \dot{H}^{-\12}(\R^n)}.
  \end{equation}
  is provided by the functions
  $f_*(x) = \Tonde{1 + \abs{x}^2}^{-(n-1)/2}$, $g_*(x) = 0$.
\end{theorem}

The geometric invariance properties of the equation~\eqref{eq:4}
suggest a way to completely characterize the class of all maximizers.

\begin{definition} \label{def:2}
  Let $\mL$ be the Lie group of transformations
  acting on \emph{solutions} of the wave equation
  and generated by:
  \begin{itemize}
  \item space-time translations:
    $u(t, x) \changeto u(t + t_0, x + x_0)$,
    with $t_0 \in \R$, $x_0 \in \R^n$;
  \item isotropic dilations:
    $u(t, x) \changeto u(\lambda t, \lambda x)$,
    with $\lambda > 0$;
  \item change of scale:
    $u(t, x) \changeto \mu u(t, x)$, with $\mu > 0$;
  \item space rotations:
      $u(t, x) \changeto u(t, R x)$, 
      with $R \in SO(n)$;
  \item phase shifts:
    $u(t, x) \changeto \Eu^{\iC \theta_+} u_+(t, x) + \Eu^{\iC \theta_-} u_-(t, x)$,
    with $\theta_+, \theta_- \in \R$
    (for the meaning of $u_+$ and $u_-$ see the next section);
  \item Lorentzian boosts:
    \begin{equation*}
      u(t, x_1, x') \changeto
       u\Tonde{\cosh(a) t + \sinh(a) x_1,
         \sinh(a) t + \cosh(a) x_1, x'},
    \end{equation*}
    with $a \in \R$.
  \end{itemize}
\end{definition}

If $u$ solves~\eqref{eq:4} and $g \in \mL$
then $v = g \cdot u$ is still a solution to~\eqref{eq:4}.
Moreover, the ratio $\norm{u}_{L^{p(n-1)}} /
\norm{(u(0), \de_t u(0))}_{\dot{H}^{\12} \times \dot{H}^{-\12}}$
is left unchanged by the action of $\mL$.

\begin{theorem}
  Let $(n, p) = (2, 6)$ or $(n, p) = (3, 4)$.
  We consider the initial data
  $f_*(x) = \Tonde{1 + \abs{x}^2}^{-(n-1)/2}$, $g_*(x) = 0$,
  and let $u_*$ be the corresponding solution
  to the wave equation~\eqref{eq:4}.
  Then the set of maximizers
  for which the equality~\eqref{eq:7} holds
  coincides with the set of initial data
  of solutions to~\eqref{eq:4} in the orbit of $u_*$
  under the action of the group $\mL$.
\end{theorem}

In order to understand how to construct maximizers,
we first present sharp proofs of the Strichartz estimates,
based on the space-time Fourier transform
in the spirit of Klainerman and Machedon's work
on bilinear estimates%
~\cite{KlaMac1996},~\cite{FosKla2000}.
We then optimize each step of the proof
by imposing conditions under which
all inequalities become equalities.
What we find are functional equations
for the Fourier transform of maximizers;
their solutions are given by particular exponential functions
with linear or quadratic exponents.

The key tool is the following well-known simple fact
about Cauchy-Schwarz's inequality for inner products.

\begin{lemma} \label{lem:0}
  Let $\ip{\cdot, \cdot}$ be a (complex) inner product on a vector
  space $V$ and let $u, v \in V$ be two non-zero vectors.
  Cauchy-Schwarz's inequality says that
  \begin{equation*}
    \abs{\ip{u, v}}^2 \le \ip{u, u} \ip{v, v};
  \end{equation*}
  moreover, equality holds if and only if $u = \alpha v$
  for some scalar $\alpha \in \C$.
\end{lemma}

\begin{remark}
  The uniqueness of maximizers modulo the transformation groups
  described in definitions~\ref{def:1} and~\ref{def:2}
  will be checked \emph{a posteriori},
  after we obtain explicit formulae for maximizers,
  and it is not used in the proof.
  While our proof relies heavily on the fact that $p$ is an even integer,
  the geometric characterization can be stated also in higher dimensions
  when $p$ is not an even integer.
  It would be interesting to prove our results
  without making use of the Fourier transform.
  For the moment, we formulate the following natural conjectures.
\end{remark}

\begin{conjecture}
  For any integer $n \ge 1$, let $p = 2 + 4/n$,
  let $f_*(x) = \exp\Tonde{-\abs{x}^2}$ 
  and $u_*(t, x) = \Eu^{-\iC t \Delta} f_*(x)$
  be the corresponding solution to the Schr\"odinger equation~\eqref{eq:2}.
  Then the set of maximizers for which the equality~\eqref{eq:6} holds
  coincides with the set of initial data
  of solutions to~\eqref{eq:2} in the orbit of $u_*$
  under the action of the group $\mG$.
\end{conjecture}

\begin{conjecture}
  For any integer $n \ge 2$,
  let $p = 2 + 4/(n-1)$ and
  let $f_*$ be the function on $\R^n$ whose Fourier transform
  is $\widehat{f_*}(\xi) = \abs{\xi}^{-1} \exp\tonde{-\abs{\xi}}$.
  Let $u_*$ be the solution to the wave equation~\eqref{eq:4}
  corresponding to the initial data
  $u_*(0) = f_*$, $\de_t u_*(0) = 0$.
  Then the set of maximizers for which the equality~\eqref{eq:7} holds
  coincides with the set of initial data
  of solutions to~\eqref{eq:4} in the orbit of $u_*$
  under the action of the group $\mL$.  
\end{conjecture}

\section{%
  Notation and preliminaries}
\label{sec:notation}

For $1 \le p < \infty$,
$L^p(\R^n)$ is the usual Lebesgue space with norm
\begin{equation*}
\norm{f}_{L^p(\R^n)} = \tonde{\int_{\R^n} \abs{f(x)}^p \d x}^{1/p}.
\end{equation*}
The homogeneous Sobolev spaces $\dot{H}^{\12}(\R^n)$ and $\dot{H}^{-\12}(\R^n)$
are defined by the norms
\begin{equation*}
  \norm{f}_{\dot{H}^{\12}(\R^n)} = \Norm{D^{\12} f}_{L^2(\R^n)}, \qquad
  \norm{f}_{\dot{H}^{-\12}(\R^n)} = \Norm{D^{-\12} f}_{L^2(\R^n)}, \qquad
\end{equation*}
where $D = \sqrt{-\Delta}$. 
In the context of the wave equation we set
\begin{equation*}
  \norm{(f, g)}_{\dot{H}^{\12}(\R^n) \times \dot{H}^{-\12}(\R^n)}
  = \tonde{\norm{f}_{\dot{H}^{\12}(\R^n)}^2
    + \norm{g}_{\dot{H}^{-\12}(\R^n)}^2}^{1/2}.
\end{equation*}

If $f(x)$ is an integrable function defined on $\R^n$,
we define its (spatial) Fourier transform by
\begin{equation*}
\fh(\xi) = \int_{\R^n} f(x) \Eu^{-\iC x \cdot \xi} \d x.
\end{equation*}
If $F(t,x)$ is an integrable function defined on $\R \times \R^n$,
we define its space-time Fourier transform by
\begin{equation*}
\Ft(\tau, \xi) =
  \int_{\R \times \R^n} F(t, x)
  \Eu^{-\iC \tonde{t \tau + x \cdot \xi}} \d t \d x.
\end{equation*}
These definitions extend in the usual way
to tempered distributions.
The Fourier transform acts like an isometry on $L^2$
and, with our definition for the Fourier transform,
Plancherel's theorem states that
\begin{equation*}
  \Norm{\fh}_{L^2(\R^n)} =
  (2 \pi)^{n/2} \norm{f}_{L^2(\R^n)}, \quad
  \Norm{\Ft}_{L^2(\R^{1 + n})} =
  (2 \pi)^{(1 + n)/2} \norm{f}_{L^2(\R^{1 + n})}.
\end{equation*}
We recall also that the Fourier transform of a pointwise product
(when it is defined)
is given by the convolution product
of the Fourier transform of each factor,
\begin{align*}
  \widehat{f g}(\xi) &=
  \1{(2 \pi)^n} \fh \conv \widehat{g} (\xi) = 
  \1{(2 \pi)^n} \int_{\R^n \times \R^n}
  \fh(\eta) \widehat{g}(\zeta)
  \ddirac{\xi - \eta - \zeta} \d\eta \d\zeta, \\
  \widetilde{F G}(\tau, \xi) &=
  \1{(2 \pi)^{n + 1}}
  \Ft \conv \Gt (\tau, \xi) = \\
  &= \1{(2 \pi)^{n + 1}}
  \int_{\R \times \R^n \times \R \times \R^n}
  \Ft(\lambda, \eta) \Gt(\mu, \zeta)
  \ddirac{\tau - \lambda - \mu \\ \xi - \eta - \zeta}
  \d\lambda \d\eta \d\mu \d\zeta.
\end{align*}
Here $\!\ddirac{\cdot}$ denotes Dirac's delta measure
concentrated in $0$,
$\int \!\ddirac{x} f(x) \d x = f(0)$.
We also denote the tensor product of two delta functions by
\begin{equation*}
  \ddirac{a \\ b} = \ddirac{a} \ddirac{b}.
\end{equation*}

\medskip
If $u(t, x)$ is the solution of the Schr\"odinger equation~\eqref{eq:2}
then its space-time Fourier transform is
\begin{equation*}
  \ut(\tau, \xi) =
  2 \pi \ddirac{\tau - \abs{\xi}^2} \fh(\xi),
\end{equation*}
where $f$ is the initial data at time $t=0$.
This shows that $\ut$ is a measure supported
on the paraboloid $\tau = \abs{\xi}^2$.
We notice,
in connection with the invariance of equation~\eqref{eq:2}
under Galilean transformations,
that the measure $\ddirac{\tau - \abs{\xi}^2}$
is invariant under the volume preserving affine change of variables
\begin{equation} \label{eq:8}
  (\tau, \xi) \changeto (\tau + 2 v \cdot \xi + \abs{v}^2, \xi + v),
\end{equation}
for any $v \in \R^n$.

\medskip
If $u(t, x)$ is the solution of the wave equation~\eqref{eq:4}
with initial data~\eqref{eq:5}, then
we can split it as $u = u_+ + u_-$, where
\begin{align*}
  u_+(t) &= \Eu^{\iC t D} D^{-\12} f_+, &
  f_+ &= \12 \tonde{D^{\12} f - \iC D^{-\12} g}, \\
  u_-(t) &= \Eu^{-\iC t D} D^{-\12} f_-, &
  f_- &= \12 \tonde{D^{\12} f + \iC D^{-\12} g}.
\end{align*}
We call $u_+$ a $(+)$-wave with data $f_+$
and $u_-$ a $(-)$-wave with data $f_-$.
Observe that, by parallelogram's law,
\begin{equation*}
  \norm{(f, g)}_{\dot{H}^{\12}(\R^n)
    \times \dot{H}^{-\12}(\R^n)}^2 = 
  2 \tonde{\norm{f_+}_{L^2(\R^n)}^2 + \norm{f_-}_{L^2(\R^n)}^2}.
\end{equation*}
The space-time Fourier transforms
of $u_+$ and $u_-$ are
\begin{equation*}
  \widetilde{u_+}(\tau, \xi) =
  2 \pi \abs{\xi}^{-\12}
  \ddirac{\tau - \abs{\xi}} \fph(\xi), \qquad
  \widetilde{u_-}(\tau, \xi) =
  2 \pi \abs{\xi}^{-\12}
  \ddirac{\tau + \abs{\xi}} \fmh(\xi).
\end{equation*}
Hence, $\widetilde{u_+}$ and $\widetilde{u_-}$
are measures supported on the null cones
$\tau = \abs{\xi}$ and $\tau = - \abs{\xi}$, respectively.
We notice also
that the measures $ \abs{\xi}^{-1} \ddirac{\tau \mp \abs{\xi}}$
are invariant under proper Lorentz transformations.
Indeed, we can write
\begin{equation*}
   \frac{\ddirac{\tau \mp \abs{\xi}}}{\abs{\xi}} = 
   2 \ddirac{\tau^2 - \abs{\xi}^2} \chic(\pm \tau > 0).
\end{equation*}

The invariance properties of these delta measures
on paraboloids and on null cones later
will help us in the computation of some convolution integrals.
Eventually we will need the following simple property of convolutions.

\begin{lemma} \label{lem:21}
  Let $A$ be a $n \times n$ invertible matrix and $b$ a vector in $\R^n$.
  Suppose $f$ is a function (or a distribution) on $\R^n$
  which is invariant under
  the linear affine change of variable $x \changeto A x + b$,
  in the sense that $f(x) = f(A x + b)$ for all $x \in \R^n$.
  Then, if the convolution $f \conv f$ is well defined,
  we have that
  \begin{equation*}
    f \conv f (x) = \1{\det A} f \conv f (A x + 2 b),
  \end{equation*}
  and, if the convolution $f \conv f \conv f$ is well defined,
  we have that
  \begin{equation*}
    f \conv f \conv f (x) = \1{(\det A)^2} f \conv f \conv f (A x + 3 b).
  \end{equation*}
\end{lemma}

For a complex number $z \in \C$,
we denote its real and imaginary parts by $\Re(z)$ and $\Im(z)$
and its complex conjugate by $\ovl{z}$.
Whenever they are mentioned,
$\log(z)$ and $\sqrt{z}$ are the branches
of the complex logarithm and of the complex square root
defined on $\C \setminus \R_-$
which extend analytically
the standard real logarithm and the standard square root
of positive real numbers.

For a vector $x = (x_1, x_2, \dots, x_n) \in \R^n$,
we adopt the prime notation to denote the vector
$x' = (x_2, \dots, x_n) \in \R^{n-1}$,
so that $x = (x_1, x')$.

If $E$ is a subset of $\R^n$
we denote its closure with respect to the usual topology
by $\ovl{E}$.

\section{%
  Schr\"odinger equation in dimension $n = 2$.}
\label{sec:case-n-2}

Consider the case $n = 2$, $p = 4$ for estimate~\eqref{eq:1}.
By Plancherel's theorem, 
$u \in L^4$ if and only if $\widetilde{u^2} \in L^2$ and
\begin{equation} \label{eq:9}
  \norm{u}_{L^4(\R^3)}^2 = \norm{u^2}_{L^2(\R^3)} = 
  (2 \pi)^{-3/2} \Norm{\widetilde{u^2}}_{L^2(\R^3)}.
\end{equation}
The Fourier transform of $u^2$ reduces to
\begin{equation} \label{eq:10}
  \widetilde{u^2}(\tau, \xi) = 
  \1{(2 \pi)^3} \ut \conv \ut (\tau, \xi) =
  \1{2 \pi} \int_{\R^2 \times \R^2}
  \fh(\eta) \fh(\zeta)
  \ddirac{\tau - \abs{\eta}^2 - \abs{\zeta}^2 \\
    \xi - \eta - \zeta}
  \d\eta \d\zeta.
\end{equation}
When $\xi = \eta + \zeta$
and $\tau = \abs{\eta}^2 + \abs{\zeta}^2$,
by parallelogram's law we have
\begin{equation*}
  2 \tau =
  \abs{\eta + \zeta}^2 + \abs{\eta - \zeta}^2 \ge \abs{\xi}^2.
\end{equation*}
It follows that $\widetilde{u^2}$ is supported
in the closure of the region
\begin{equation*}
  \mP_2 = \graffe{(\tau, \xi) \in \R \times \R^2:
    2 \tau > \abs{\xi}^2}.
\end{equation*}
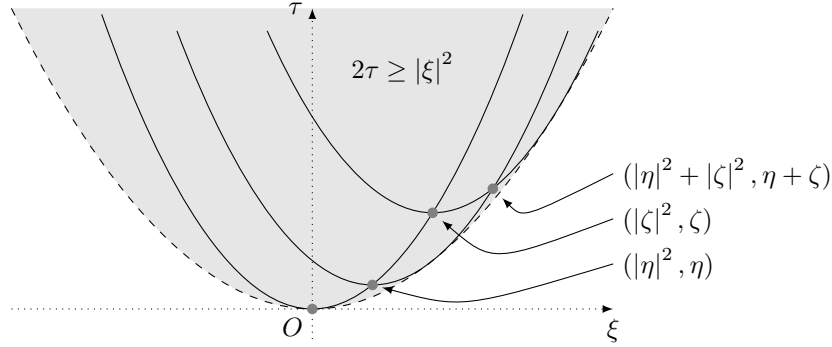
\begin{figure}
  \centering
  \begin{tikzpicture}[scale=2, >=latex]
    \fill[black!10]
    (-2, 2) parabola[bend at end] (0,0) parabola (2, 2) -- cycle;
    \draw[dashed]
    (-2, 2) parabola[bend at end] (0,0) parabola (2, 2);
    \draw[->, thin, dotted] (-2, 0) -- (2, 0) node[below] {$\xi$};
    \draw[->, thin, dotted] (0, -0.2) -- (0, 2) node[left] {$\tau$};
    \draw
    (-1.4, 1.96) parabola[bend at end] (0,0) parabola (1.4, 1.96);
    \draw[shift={(0.8,0.64)}]
    (-1.1, 1.21) parabola[bend at end] (0,0) parabola (1.1, 1.21);
    \draw[shift={(0.4,0.16)}]
    (-1.3, 1.69) parabola[bend at end] (0,0) parabola (1.3, 1.69);
    \fill[black!50, text=black]
    (0, 0) circle (1pt) node[below left] {$O$}
    (0.4, 0.16) node (a) {} circle (1pt)
    (0.8, 0.64) node (b) {} circle (1pt)
    (1.2, 0.80) node (c) {} circle (1pt);
    \path (0.6, 1.6) node {$2\tau \ge \abs{\xi}^2$};
    \draw[->]
    (2, 0.3) node[right] {$(\abs{\eta}^2, \eta)$}
    .. controls (1, 0) ..  (a);
    \draw[->]
    (2, 0.6) node[right] {$(\abs{\zeta}^2, \zeta)$}
    .. controls (1.2, 0.3) ..  (b);
    \draw[->]
    (2, 0.9) node[right] {$(\abs{\eta}^2 + \abs{\zeta}^2, \eta + \zeta)$}
    .. controls (1.4, 0.6) ..  (c);
  \end{tikzpicture}
  \caption{The region $\mP_2$}
  \label{fig:5}
\end{figure}
For each choice of $(\tau, \xi) \in \mP_2$,
we denote by $\ip{\cdot, \cdot}_{(\tau, \xi)}$
the $L^2$ inner product associated with the measure
\begin{equation} \label{eq:11}
  \MU_{(\tau, \xi)} =
  \ddirac{\tau - \abs{\eta}^2 - \abs{\zeta}^2 \\ \xi - \eta - \zeta}
  \d\eta \d\zeta;
\end{equation}
and by $\norm{\cdot}_{(\tau, \xi)}$ the corresponding norm;
more precisely, we set
\begin{align*}
  & \ip{F, G}_{(\tau, \xi)} =
  \int_{\R^2 \times \R^2} F(\eta, \zeta) \ovl{G(\eta, \zeta)}
  \ddirac{\tau - \abs{\eta}^2 - \abs{\zeta}^2 \\
    \xi - \eta - \zeta}
  \d\eta \d\zeta, \\
  & \norm{F}_{(\tau, \xi)} =
  \tonde{\int_{\R^2 \times \R^2} \abs{F(\eta, \zeta)}^2 
    \ddirac{\tau - \abs{\eta}^2 - \abs{\zeta}^2 \\
      \xi - \eta - \zeta}
    \d\eta \d\zeta}^{1/2}.
\end{align*}

\begin{remark}
  The measure $\MU_{(\tau, \xi)}$ defined in~\eqref{eq:11}
  is the \emph{pull-back} of the Dirac's delta on $\R \times \R^2$
  by the function $\Phi_{(\tau, \xi)}: \R^2 \times \R^2 \to \R \times \R^2$
  given by
  \begin{equation*}
    \Phi_{(\tau, \xi)}(\eta, \zeta) =
    \Tonde{\tau - \abs{\eta}^2 - \abs{\zeta}^2, \xi - \eta - \zeta}.
  \end{equation*}
  This pull back is well defined as long as
  the differential of $\Phi_{(\tau, \xi)}$ is surjective
  in correspondence of the points where $\Phi_{(\tau, \xi)}$ vanishes
  (we refer to~\cite[Theorem 6.1.2 and Example 6.1.3]{Hor1990}
  for more details about pull-backs of distributions).
  A quick computation shows that the differential of $\Phi_{(\tau, \xi)}$
  is surjective at a point $(\eta, \zeta)$
  if and only if $\eta \ne \zeta$.
  On the other hand,
  if $\Phi_{(\tau, \xi)}(\eta, \eta) = 0$
  we must have
  \begin{equation*}
    2\tau = 2(\abs{\eta}^2 + \abs{\eta}^2) = \abs{\eta + \eta}^2 = \abs{\xi}^2.
  \end{equation*}  
  This tells us that $\MU_{(\tau, \xi)}$ is not well defined
  on the boundary of $\mP_2$, when $2 \tau = \abs{\xi}^2$.
  However, we can safely ignore the problems at this boundary
  and observe instead that
  for any locally integrable function $F(\eta, \zeta)$
  defined on $\R^2 \times \R^2$
  the integral $G(\tau, \xi) = \int F \d\MU_{(\tau, \xi)}$
  defines a locally integrable function on $\R \times \R^2$.
  Indeed, if $K$ is a compact set in $\R \times \R^2$, we have
  \begin{multline*}
    \iint_K \abs{G(\tau, \xi)} \d\tau \d\xi \le
    \iint_{(\tau, \xi) \in K} \iint \abs{F(\eta, \zeta)} 
    \ddirac{\tau - \abs{\eta}^2 - \abs{\zeta}^2 \\ \xi - \eta - \zeta}
    \d\eta \d\zeta \d\tau \d\xi = \\
    = \iint_{(\abs{\eta}^2 + \abs{\zeta}^2, \eta + \zeta) \in K}
    \abs{F(\eta, \zeta)} \d\eta \d\zeta,
  \end{multline*}
  and the set 
  $\Graffe{(\eta, \zeta): (\abs{\eta}^2 + \abs{\zeta}^2, \eta + \zeta) \in K}$
  is a compact set in $\R^2 \times \R^2$.
\end{remark}

We can now write~\eqref{eq:10} as
\begin{equation*}
  \widetilde{u^2}(\tau, \xi) = \1{2 \pi}
  \ip{\fh \otimes \fh, 1 \otimes 1}_{(\tau, \xi)},
\end{equation*}
where the tensor product is defined by
$(f \otimes g) (\eta, \zeta) = f(\eta) g(\zeta)$.
By Cauchy-Schwarz's inequality we obtain that
\begin{equation} \label{eq:12}
  \Abs{\widetilde{u^2}(\tau, \xi)} \le \1{2 \pi}
  \Norm{\fh \otimes \fh}_{(\tau, \xi)}
  \Norm{1 \otimes 1}_{(\tau, \xi)}.
\end{equation}
Hence,
\begin{equation} \label{eq:13}
  \Norm{\widetilde{u^2}}_{L^2(\R^3)} \le \1{2 \pi}
  \TOnde{\sup_{(\tau, \xi) \in \mP_2}
    \Norm{1 \otimes 1}_{(\tau, \xi)}}
  \tonde{\int_{\mP_2}
    \Norm{\fh \otimes \fh}_{(\tau, \xi)}^2
    \d\tau \d\xi}^{1/2}.
\end{equation}
The next lemma shows that
not only $\Norm{1 \otimes 1}_{(\tau, \xi)}$ is uniformly bounded
with respect to $(\tau, \xi)$,
but that it is actually constant
on the support of $\widetilde{u^2}$.

\begin{lemma} \label{lem:1}
  For each $(\tau, \xi) \in \mP_2$ we have 
  $\norm{1 \otimes 1}_{(\tau, \xi)} = \sqrt{\pi/2}$.
\end{lemma}

\begin{proof}
  The quantity
  \begin{equation*}
    I(\tau, \xi) = \norm{1 \otimes 1}_{(\tau, \xi)}^2 =
    \int_{\R^2} \ddirac{\tau - \abs{\xi - \eta}^2 - \abs{\eta}^2} \d\eta
  \end{equation*}
  is just the convolution of the measure $\ddirac{\tau - \abs{\xi}^2}$
  with itself.
  The invariance of this measure
  with respect to the transformation~\eqref{eq:8}
  together with lemma~\ref{lem:21}
  imply that
  \begin{equation*}
    I(\tau, \xi) = I(\tau + 2 v \cdot \xi + \abs{v}^2, \xi + v),
  \end{equation*}
  for any $v \in \R^2$.
  If we take $v = -\xi/2$ we obtain that $I(\tau, \xi) = I(\tau^*, 0)$,
  where $\tau^* = \tau - \abs{\xi}^2/2$.
  Moreover, it is evident from its definition that,
  by homogeneity,
  $I$ is invariant under parabolic dilations,
  $I(\tau, \xi) = I(\lambda^2 \tau, \lambda \xi)$.
  Hence, when $\tau^* > 0$ we have
  \begin{multline*}
    I(\tau, \xi) = I(\tau^*, 0) = I(1, 0) = 
    \int_{\R^2} \ddirac{1 - 2\abs{\eta}^2} \d\eta = \\
    = 2 \pi \int_0^\infty \ddirac{1 - 2 r^2} r \d r =
    \pi \int_0^\infty \ddirac{1 - 2 s} \d s =
    \frac\pi2.
  \end{multline*}
\end{proof}

We also have
\begin{multline} \label{eq:14}
  \int_{\mP_2}
  \Norm{\fh \otimes \fh}_{(\tau, \xi)}^2
  \d\tau \d\xi = \\
  = \int_{\R^2 \times \R^2}
  \abs{\fh(\eta) \fh(\zeta)}^2
  \int_{\mP_2}
  \ddirac{\tau - \abs{\eta}^2 - \abs{\zeta}^2 \\
    \xi - \eta - \zeta}
  \d\tau \d\xi \d\eta \d\zeta = \\
  = \Norm{\fh \otimes \fh}_{L^2(\R^2 \times \R^2)}^2
  = \Norm{\fh}_{L^2(\R^2)}^4
  = (2 \pi)^4 \norm{f}_{L^2(\R^2)}^4.
\end{multline}
It follows from~\eqref{eq:9},~\eqref{eq:13},~\eqref{eq:14}
and lemma~\ref{lem:1} that
\begin{equation} \label{eq:15}
  \norm{u}_{L^4(\R^3)} \le \1{\sqrt{2}} \norm{f}_{L^2(\R^2)}.
\end{equation}
This proves that for $n = 2$
the best constant $S(2)$ in~\eqref{eq:1}
is no larger that $1/\sqrt{2}$.

\begin{remark}
  We observe that in the above computations
  the only place where we have used an inequality
  instead of an equality is in~\eqref{eq:13}
  as a consequence of the Cauchy-Schwarz inequality~\eqref{eq:12}.
  If we can find a function $f$
  for which we have equality in~\eqref{eq:12}
  for all $(\tau, \xi) \in \mP_2$
  then there will be equality also in~\eqref{eq:15}.
  This will show that $f$ is a maximizer for the estimate
  and that $S(2) = 1/\sqrt2$.
\end{remark}

We have equality in~\eqref{eq:13}
if there is equality in~\eqref{eq:12}
for almost all $(\tau, \xi) \in \mP_2$.
By lemma~\ref{lem:0}, this happens
if there exists a scalar function $F: \mP_2 \to \C$ such that
\begin{equation*}
  \Tonde{\fh \otimes \fh}(\eta, \zeta) =
  F(\tau, \xi) \tonde{1 \otimes 1}(\eta,\zeta),
\end{equation*}
for almost all $(\eta, \zeta)$
(with respect to the measure~\eqref{eq:11})
in the support of the measure~\eqref{eq:11}
and for almost all $(\tau, \xi)$ in $\mP_2$
(with respect to the Lebesgue measure on~$\R \times \R^2$).
This means that we are looking for functions $f$ and $F$ such that
\begin{equation} \label{eq:16}
  \fh(\eta) \fh(\zeta) =
  F\tonde{\abs{\eta}^2 + \abs{\zeta}^2, \eta + \zeta},
\end{equation}
for almost all $(\eta, \zeta) \in \R^2 \times \R^2$.
An example of such functions is given by the pair
$\fh(\xi) = \Eu^{-\abs{\xi}^2}$, $F(\tau, \xi) = \Eu^{-\tau}$.

If $f$ is a maximizer, $\fh$ must solve the equation~\eqref{eq:16}
and it follows from proposition~\ref{pro:3} that
\begin{equation} \label{eq:17}
  \fh(\xi) = \exp\Tonde{\Ah \abs{\xi}^2 + \bh \cdot \xi + \Ch},
  \qquad \xi \in \R^2,
\end{equation}
for some constants 
$\Ah \in \C$,
$\bh = \Tonde{\bh_1, \bh_2} \in \C^2$,
$\Ch \in \C$,
with $\Re(\Ah) < 0$ in order to have $f \in L^2(\R^2)$.
The inverse Fourier transform of~\eqref{eq:17}
is again a function of the same class
\begin{equation} \label{eq:18}
  f(x) = \exp\Tonde{A \abs{x}^2 + b \cdot x + C}, \qquad
  x \in \R^2,
\end{equation}
where the relations between
the parameters $A \in \C$, $b \in \C^2$, $C \in \C$
and the parameters $\Ah$, $\bh$, $\Ch$
are given by
\begin{equation*}
  A = \1{4\Ah}, \qquad
  b = -\frac{\iC \bh}{4 \Ah}, \qquad
  C = \Ch - \frac{\bh_1^2 + \bh_2^2}{4 \Ah} - \log(-4 \pi \Ah).
\end{equation*}

The class of initial data of the form~\eqref{eq:18}
is invariant under the action of the group~$\mG$
described in definition~\ref{def:1}.
The coefficients change according to the following rules:
  \begin{itemize}
  \item space-time translations:
    $(\Ah, \bh, \Ch) \changeto
    (\Ah + \iC t_0, \bh + \iC x_0, \Ch)$;
  \item parabolic dilations:
    $(\Ah, \bh, \Ch) \changeto
    (\Ah / \lambda^2, \bh / \lambda, \Ch - n \log\lambda)$;
  \item change of scale:
    $(\Ah, \bh, \Ch) \changeto (\Ah, \bh, \Ch + \log\mu)$;
  \item space rotations:
    $(\Ah, \bh, \Ch) \changeto (\Ah, R \bh, \Ch)$;
  \item phase shifts:
    $(\Ah, \bh, \Ch) \changeto (\Ah, \bh, \Ch + \iC \theta)$;
  \item Galilean transformations:
    \begin{equation*}
    (\Ah, \bh, \Ch) \changeto
    \TOnde{\Ah, \bh + A v,
      \Ch - \frac{\Ah}4 \abs{v}^2 - \frac{\bh}2 \cdot v}.
  \end{equation*}
\end{itemize}

Hence, after a translation and a phase shift
we can make all coefficients real;
by a Galilean transformation we can make $\bh = 0$;
then, by a parabolic dilation we can have $\Ah = -1/4$;
finally a change of scale gives $\Ch = \log(-\pi)$.
This would correspond to the case $A = -1$, $b = 0$, $C = 0$,
which is the function $f_*(x) = \Eu^{-\abs{x}^2}$.
Thus, we have proved that any maximizer is connected to $f_*$
by the action of $\mG$.

\section{%
  Schr\"odinger equation in dimension $n = 1$.}
\label{sec:case-n-1}

Consider the case $n = 1$, $p = 6$ for estimate~\eqref{eq:1}.
This is the case that was considered in \cite{Kun2003}.
By Plancherel's theorem,
$u \in L^6$ if and only if $\widetilde{u^3} \in L^2$ and
\begin{equation} \label{eq:19}
  \norm{u}_{L^6(\R^2)}^3 = \norm{u^3}_{L^2(\R^2)} = 
  (2 \pi)^{-1} \Norm{\widetilde{u^3}}_{L^2(\R^2)}.
\end{equation}
The Fourier transform of $u^3$ reduces to
\begin{multline} \label{eq:20}
  \widetilde{u^3}(\tau, \xi) = \1{(2 \pi)^4}
  \ut \conv \ut \conv \ut (\tau, \xi) = \\
  = \1{2 \pi} \int_{\R \times \R \times \R}
  \fh(\eta_1) \fh(\eta_2) \fh(\eta_3)
  \ddirac{\tau - \eta_1^2 - \eta_2^2 - \eta_3^2 \\
    \xi - \eta_1 - \eta_2 - \eta_3}
  \d\eta_1 \d\eta_2 \d\eta_3 = \\
  = \1{2 \pi} \int_{\R^3}
  \fh(\eta_1) \fh(\eta_2) \fh(\eta_3)
  \ddirac{\tau - \abs{\eta}^2 \\ \xi - (1, 1, 1) \cdot \eta}
  \d\eta,
\end{multline}
where now $\eta = \tonde{\eta_1, \eta_2, \eta_3}$.
When $\xi = (1, 1, 1) \cdot \eta$ and $\tau = \abs{\eta}^2$,
we have $3 \tau \ge \xi^2$.
It follows that $\widetilde{u^3}$ is supported
in the closure of the region
\begin{equation*}
  \mP_1 = \graffe{(\tau, \xi) \in \R \times \R: 3 \tau > \xi^2}.
\end{equation*}
For each choice of $(\tau, \xi) \in \mP_1$,
we denote by $\ip{\cdot, \cdot}_{(\tau, \xi)}$
the $L^2$ inner product associated with the measure
\begin{equation} \label{eq:21}
  \ddirac{\tau - \abs{\eta}^2 \\
    \xi - (1, 1, 1) \cdot \eta}
  \d\eta;
\end{equation}
and by $\norm{\cdot}_{(\tau, \xi)}$ the corresponding norm.
We can then write~\eqref{eq:20} as
\begin{equation*}
  \widetilde{u^3}(\tau, \xi) = \1{2 \pi}
  \ip{\fh \otimes \fh \otimes \fh,
    1 \otimes 1 \otimes 1}_{(\tau, \xi)},
\end{equation*}
where the tensor product is defined by
$(f \otimes g \otimes h) (\eta) = f(\eta_1) g(\eta_2) h(\eta_3)$.
By Cauchy-Schwarz's inequality we obtain that
\begin{equation} \label{eq:22}
  \Abs{\widetilde{u^3}(\tau, \xi)} \le \1{2 \pi}
  \Norm{\fh \otimes \fh \otimes \fh}_{(\tau, \xi)}
  \Norm{1 \otimes 1 \otimes 1}_{(\tau, \xi)}.
\end{equation}
Hence,
\begin{equation} \label{eq:23}
  \Norm{\widetilde{u^3}}_{L^2(\R^2)} \le \1{2 \pi}
  \TOnde{\sup_{(\tau, \xi) \in \mP_1}
    \Norm{1 \otimes 1 \otimes 1}_{(\tau, \xi)}}
  \tonde{\int_{\mP_1} \Norm{\fh \otimes \fh \otimes \fh}_{(\tau, \xi)}^2
    \d\tau \d\xi}^{1/2}.
\end{equation}
The next lemma shows that
not only $\Norm{1 \otimes 1 \otimes 1}_{(\tau, \xi)}$ is bounded,
but that it is actually constant
on the support of $\widetilde{u^3}$.

\begin{lemma} \label{lem:2}
  For each $(\tau, \xi) \in \mP_1$ we have 
  $\norm{1 \otimes 1 \otimes 1}_{(\tau, \xi)}
  = \sqrt{\pi / \sqrt3}$.
\end{lemma}

\begin{proof}
  The quantity
  \begin{equation*}
    I(\tau, \xi) = \norm{1 \otimes 1 \otimes 1}_{(\tau, \xi)}^2 =
    \int_{\R^3}
    \ddirac{\tau - \abs{\eta}^2 \\ \xi - (1, 1, 1) \cdot \eta}
    \d\eta
  \end{equation*}
  is just an iterated convolution of the measure $\ddirac{\tau - \xi^2}$
  with itself twice.
  {}From lemma~\ref{lem:21} and the invariance of this measure
  with respect to the transformation~\eqref{eq:8},
  it follows that $I(\tau, \xi) = I(\tau^*, 0)$
  where $\tau^* = \tau - \xi^2/3$.
  Moreover, by homogeneity $I$ is invariant under parabolic dilations,
  $I(\lambda^2 \tau, \lambda \xi) = I(\tau, \xi)$.
  Hence, when $\tau^* > 0$ we have
  \begin{multline*}
    I(\tau, \xi) = I(\tau^*, 0) = I(1, 0) = 
    \int_{\R^3} \ddirac{1 - \abs{\eta}^2 \\ -(1, 1, 1) \cdot \eta} \d\eta =
    \int_{\R^3}
    \ddirac{1 - \abs{\eta}^2 \\ \abs{(1, 1, 1)} \eta_1} \d\eta = \\
    = \1{\sqrt{3}} \int_{\R^2} \ddirac{1 - \abs{\zeta}^2} \d\zeta =
    \frac{2\pi}{\sqrt{3}} \int_0^\infty \ddirac{1 - r^2} r \d r = 
    \frac\pi{\sqrt{3}}.
  \end{multline*}
\end{proof}

We also have
\begin{multline} \label{eq:24}
  \int_{\mP_1} \Norm{\fh \otimes \fh \otimes \fh}_{(\tau, \xi)}^2
  \d\tau \d\xi = \\
  = \int_{\R^3}
  \abs{\fh(\eta_1) \fh(\eta_2) \fh(\eta_3)}^2
  \int_{\mP_1} \ddirac{\tau - \abs{\eta}^2 \\ \xi - (1, 1, 1) \cdot \eta}
  \d\tau \d\xi \d\eta = \\
  = \Norm{\fh \otimes \fh \otimes \fh}_{L^2(\R^3)}^2
  = \Norm{\fh}_{L^2(\R)}^6
  = (2 \pi)^3 \norm{f}_{L^2(\R)}^6.
\end{multline}
It follows from~\eqref{eq:19},~\eqref{eq:23},~\eqref{eq:24}
and lemma~\ref{lem:2} that
\begin{equation} \label{eq:25}
  \norm{u}_{L^6(\R^2)} \le {12}^{-1/12} \norm{f}_{L^2(\R)}.
\end{equation}
This proves that for $n=1$ the best constant $S(1)$ in~\eqref{eq:1}
is no larger that ${12}^{-1/12}$.

As before, we observe that
if we could find a function $f$
for which we have equality in~\eqref{eq:22}
for all $(\tau, \xi) \in \mP_1$
then we would have equality in~\eqref{eq:25}
and we would have found a maximizer for the estimate.
We have equality in the Cauchy-Schwarz inequality~\eqref{eq:22}
for (almost) all $(\tau, \xi) \in \mP_1$
if there exists a scalar function $F: \mP_1 \to \C$ such that
\begin{equation*}
  \Tonde{\fh \otimes \fh \otimes \fh}(\eta) =
  F(\tau, \xi) \tonde{1 \otimes 1 \otimes 1}(\eta),
\end{equation*}
for (almost) all $\eta$
in the support of the measure~\eqref{eq:21}.
This means that we are looking for functions $f$ and $F$ such that
\begin{equation} \label{eq:26}
  \fh(\eta_1) \fh(\eta_2) \fh(\eta_3) =
  F\Tonde{\eta_1^2 + \eta_2^2 + \eta_3^2, \eta_1 + \eta_2 + \eta_3},
\end{equation}
for (almost) all $\eta \in \R^3$.
Again, an example of such functions is given by the pair
$\fh(\xi) = \Eu^{-\xi^2}$, $F(\tau, \xi) = \Eu^{-\tau}$.

If $f$ is a maximizer, $\fh$ must solve the equation~\eqref{eq:26}
and it follows from proposition~\ref{pro:1} that
\begin{equation} \label{eq:27}
  \fh(\xi) = \exp\Tonde{\Ah \xi^2 + \Bh \xi + \Ch},
  \qquad \xi \in \R,
\end{equation}
for some complex constants $\Ah$, $\Bh$, $\Ch$,
with $\Re(\Ah) < 0$ in order to have $f \in L^2(\R)$.
The inverse Fourier transform of~\eqref{eq:27}
is again a function of the same class
\begin{equation} \label{eq:28}
  f(x) = \exp\Tonde{A x^2 + B x + C}, \qquad x \in \R,
\end{equation}
where the relations between
the parameters $A$, $B$, $C$ and the parameters $\Ah$, $\Bh$, $\Ch$
are given by
\begin{equation*}
  A = \1{4 \Ah}, \qquad
  B = -\frac{\iC \Bh}{4 \Ah}, \qquad
  C = \Ch - \frac{\Bh^2}{4 \Ah} - \12\log(-4 \pi \Ah).
\end{equation*}

As we have seen at the end of section~\ref{sec:case-n-2},
the class of initial data of the form~\eqref{eq:28}
is invariant under the action of the group~$\mG$
and any maximizer is connected to
the function $f_*(x) = \Eu^{-x^2}$
by the action of $\mG$.

\section{%
  Wave equation in dimension $n = 3$.}
\label{sec:case-w-n-3}

Consider the case $n = 3$, $p = 4$ for estimate~\eqref{eq:3}.
We have
\begin{equation} \label{eq:29}
  \widetilde{u_+ u_+}(\tau, \xi) =
  \1{(2 \pi)^2} \int_{\R^3 \times \R^3}
  \frac{\fph(\eta) \fph(\zeta)}
  {\abs{\eta}^{\12} \abs{\zeta}^{\12}}
  \ddirac{\tau - \abs{\eta} - \abs{\zeta} \\
    \xi - \eta - \zeta}
  \d\eta \d\zeta.
\end{equation}
In particular, it follows that
the term $\widetilde{u_+ u_+}$ is supported
in the closure of the region
\begin{equation*}
  \mC_{++} = \graffe{(\tau, \xi) \in \R \times \R^3: \tau > \abs{\xi}}.
\end{equation*}
Similarly, the term $\widetilde{u_+ u_-}$ is supported
in the closure of the region
\begin{equation*}
  \mC_{+-} =
  \graffe{(\tau, \xi) \in \R \times \R^3:
    \abs{\tau} < \abs{\xi}};
\end{equation*}
and the term $\widetilde{u_- u_-}$ is supported
in the closure of the region
\begin{equation*}
  \mC_{--} =
  \graffe{(\tau, \xi) \in \R \times \R^3: \tau < -\abs{\xi}}.
\end{equation*}

We remark that formulae like~\eqref{eq:29} are the starting point
for the bilinear estimates studied in~\cite{FosKla2000}.

We first prove the estimate for $u_+$.
By Plancherel's theorem we have
\begin{equation} \label{eq:30}
  \norm{u_+}_{L^4(\R^4)}^2 = \norm{u_+^2}_{L^2(\R^4)} = 
  (2 \pi)^{-2} \Norm{\widetilde{u_+^2}}_{L^2(\R^4)}.
\end{equation}
For each choice of $(\tau, \xi) \in \mC_{++}$,
we denote by $\ip{\cdot, \cdot}_{(\tau, \xi)}$
the $L^2$ inner product associated with the measure
\begin{equation} \label{eq:31}
  \MU_{(\tau, \xi)} =
  \ddirac{\tau - \abs{\eta} - \abs{\zeta} \\ \xi - \eta - \zeta}
  \d\eta \d\zeta
\end{equation}
and by $\norm{\cdot}_{(\tau, \xi)}$ the corresponding norm.

\begin{remark} \label{rem:7}
  The measure $\MU_{(\tau, \xi)}$ defined in~\eqref{eq:31}
  is the \emph{pull-back} of the Dirac's delta on $\R^{1 + 3}$
  by the function $\Phi_{(\tau, \xi)}:
  (\R^3 \setminus 0) \times (\R^3 \setminus 0) \to \R \times \R^3$
  given by
  \begin{equation*}
    \Phi_{(\tau, \xi)}(\eta, \zeta) =
    \Tonde{\tau - \abs{\eta} - \abs{\zeta}, \xi - \eta - \zeta}.
  \end{equation*}
  A quick computation shows that the differential of $\Phi_{(\tau, \xi)}$
  is surjective at a point~$(\eta, \zeta)$
  if and only if $\eta / \abs{\eta} \ne \zeta / \abs{\zeta}$.
  On the other hand,
  if $\eta / \abs{\eta} = \zeta / \abs{\zeta}$
  and~$\Phi_{(\tau, \xi)}(\eta, \zeta) = 0$
  we must have
  $\tau = \abs{\eta} + \abs{\zeta} = \abs{\eta + \zeta} = \abs{\xi}$.  
  This tells us that~$\MU_{(\tau, \xi)}$ is not well defined
  on the boundary of $\mC_{++}$, when $\tau = \abs{\xi}$.
  However, we can safely ignore the problems at this boundary
  and observe instead that
  for any locally integrable function $F(\eta, \zeta)$
  defined on $\R^3 \times \R^3$
  the integral $G(\tau, \xi) = \int F \d\MU_{(\tau, \xi)}$
  defines a locally integrable function on $\R \times \R^3$.
  Indeed, if $K$ is a compact set in $\R \times \R^3$, we have
  \begin{multline*}
    \iint_K \abs{G(\tau, \xi)} \d\tau \d\xi \le
    \iint_{(\tau, \xi) \in K} \iint \abs{F(\eta, \zeta)}
    \ddirac{\tau - \abs{\eta} - \abs{\zeta} \\ \xi - \eta - \zeta}
    \d\eta \d\zeta \d\tau \d\xi = \\
    = \iint_{(\abs{\eta} + \abs{\zeta}, \eta + \zeta) \in K}
    \abs{F(\eta, \zeta)} \d\eta \d\zeta,
  \end{multline*}
  and the set 
  $\Graffe{(\eta, \zeta): (\abs{\eta} + \abs{\zeta}, \eta + \zeta) \in K}$
  is a compact set in $\R^3 \times \R^3$.
\end{remark}

We can now write~\eqref{eq:29} as
\begin{equation} \label{eq:32}
  \widetilde{u_+^2}(\tau, \xi) = \1{(2 \pi)^2}
  \ip{\fph \otimes \fph,
    \abs{\cdot}^{-\12} \otimes \abs{\cdot}^{-\12}}_{(\tau, \xi)}.
\end{equation}
The quantity we want to compute this time is
$\norm{\abs{\cdot}^{-\12} \otimes
  \abs{\cdot}^{-\12}}_{(\tau, \xi)}$.

\begin{lemma} \label{lem:3}
  For each $(\tau, \xi) \in \mC_{++}$ we have 
  $\norm{\abs{\cdot}^{-\12} \otimes
    \abs{\cdot}^{-\12}}_{(\tau, \xi)} = (2 \pi)^{1/2}$.
\end{lemma}

\begin{proof}
  The quantity
  \begin{equation*}
    I(\tau, \xi) =
    \norm{\abs{\cdot}^{-\12} \otimes \abs{\cdot}^{-\12}}_{(\tau, \xi)}^2 =
    \int_{\R^3}
    \frac{\ddirac{\tau - \abs{\xi - \eta} - \abs{\eta}}}
    {\abs{\xi - \eta} \abs{\eta}} \d\eta
  \end{equation*}
  is just the convolution
  of the measure $\abs{\xi}^{-1} \ddirac{\tau - \abs{\xi}}$
  with itself.
  If $\tau > \abs{\xi}$, from the invariance of this measure
  with respect to proper Lorentz transformations
  and from the fact that it is always possible
  to find a proper Lorentz transformation which
  takes $(\tau, \xi)$ to the point $(\tau^*, 0)$
  where $\tau^* = (\tau^2 - \abs{\xi}^2)^{1/2}$,
  it follows that $I(\tau, \xi) = I(\tau^*, 0)$.
  Moreover, it is evident from its definition that
  by homogeneity~$I$ is invariant under isotropic dilations,
  $I(\lambda \tau, \lambda \xi) = I(\tau, \xi)$.
  Hence, when $\tau > \abs{\xi}$,
  \begin{equation*}
    I(\tau, \xi) = I(1, 0) = 
    \int_{\R^3} \frac{\ddirac{1 - 2 \abs{\xi}}}{\abs{\xi}^2} \d\xi =
    4 \pi \int_0^\infty \ddirac{1 - 2 r} \d r = 2 \pi.
  \end{equation*}
\end{proof}

Cauchy-Schwarz's inequality applied to~\eqref{eq:32}
together with lemma~\ref{lem:3} give
\begin{multline} \label{eq:33}
  \Norm{\widetilde{u_+^2}}_{L^2(\R^4)}^2 \le
  \1{(2 \pi)^3} \int_{\mC_{++}}
  \norm{\fph \otimes \fph}_{(\tau, \xi)}^2 \d\tau \d \xi = \\
  = \1{(2 \pi)^3}
  \norm{\fph \otimes \fph}_{L^2(\R^3 \times \R^3)}^2
  = \1{(2 \pi)^3} \norm{\fph}_{L^2(\R^3)}^4 =
  (2 \pi)^3 \norm{f_+}_{L^2(\R^3)}^4.
\end{multline}
Hence, combining~\eqref{eq:30} and~\eqref{eq:33} we obtain
\begin{equation} \label{eq:34}
  \Norm{u_+}_{L^4(\R^4)} \le (2 \pi)^{-1/4} \norm{f_+}_{L^2(\R^3)}.
\end{equation}
This time,
equality holds if there exist a function $F: \mC_{++} \to \C$
such that 
\begin{equation*}
  \Tonde{\fph \otimes \fph}(\eta, \zeta) =
  F(\tau, \xi) \abs{\eta}^{-\12} \abs{\zeta}^{-\12},
\end{equation*}
for all $(\eta, \zeta)$
in the support of the measure~\eqref{eq:31}.
This means that
\begin{equation*}
  \abs{\eta}^{\12} \fph(\eta) \abs{\zeta}^{\12} \fph (\zeta)
  = F\Tonde{\abs{\eta} + \abs{\zeta}, \eta + \zeta},
\end{equation*}
for almost all $\eta, \zeta \in \R^3$.
An example of such functions is given by the pair
\begin{equation*}
  \fph(\xi) = \abs{\xi}^{-\12} \Eu^{-\abs{\xi}}, \qquad
  F(\tau, \xi) = \Eu^{-\tau}.
\end{equation*}
It follows from proposition~\ref{pro:4} that
any maximizer for the estimate~\eqref{eq:34}
is a function whose Fourier transform has the form
\begin{equation} \label{eq:35}
  \fh(\xi) =
  \abs{\xi}^{-\12} \exp\TOnde{A \abs{\xi} + b \cdot \xi + C},
\end{equation}
with $A, C \in \C$, $b \in \C^3$, $\Im(C) \in [0, 2\pi[$
and $\abs{\Re(b)} < -\Re(A)$ (in order to have $f \in L^2(\R^3)$).
In the next lemma we compute an explicit expression
for homogeneous waves with data of the form~\eqref{eq:35}.

\begin{lemma} \label{lem:6}
  Let $u$ be the $(+)$-wave corresponding
  to an $L^2$ data of the form~\eqref{eq:35},
  \begin{equation} \label{eq:36}
    u(t, x) = \1{(2\pi)^3} \int_{\R^3}
    \exp\TOnde{(A + \iC t) \abs{\xi} + (b + \iC x) \cdot \xi + C}
    \frac{\d\xi}{\abs{\xi}}.
  \end{equation}
  Then we have the explicit formula
  \begin{multline} \label{eq:37}
    2 \pi^2 \Eu^{-\iC \Im(C)} u\Tonde{t - \Im(A), x - \Im(b)} = \\
    = \frac{\Eu^{\Re(C)}}
    {(\Re(A))^2 - \abs{\Re(b)}^2 + \abs{x}^2 - t^2 +
      2 \iC \tonde{\Re(A) t - \Re(b) \cdot x}}.
  \end{multline}
\end{lemma}

\begin{proof}
  The integral
  \begin{equation*}
    F(t, x) = \int_{\R^3} \exp\Tonde{t \abs{\xi} + x \cdot \xi}
    \frac{\d\xi}{\abs{\xi}}
  \end{equation*}
  is well defined for $t \in \C$ and $x \in \C^3$
  when $\Re(t) < -\abs{\Re(x)}$.
  For $t \in \R$ and $x \in \R^3$ with $t < -\abs{x}$,
  using polar coordinates,
  $r = \abs{\xi}$ and $u = (x/\abs{x}) \cdot (\xi/\abs{\xi})$,
  we find
  \begin{equation*}
    F(t, x) = 2 \pi \int_{-1}^1 \int_0^\infty
    \exp\Tonde{(t + \abs{x} u) r} r \d r \d u =
    \int_{-1}^1 \frac{2 \pi \d u}{\tonde{t + \abs{x} u}^2}
    = \frac{4 \pi}{t^2 - x_1^2 - x_2^2 - x_3^2}.
  \end{equation*}
  By analytic continuation this formula remains valid
  for complex $t$ and $x$ with $\Re(t) < -\abs{\Re(x)}$.
  Formula~\eqref{eq:37} follows from the identity  
  \begin{equation*}
    \Eu^{-\iC \Im(C)} u\Tonde{t - \Im(A), x - \Im(b)}
    = \frac{\Eu^{\Re(C)}}{(2 \pi)^3} 
    F\Tonde{\Re(A) + \iC t, \Re(b) + \iC x}.
  \end{equation*}
\end{proof}

\begin{remark} \label{rem:2}
  If $u$ is the $(+)$-wave corresponding
  to an $L^2$ data of the form~\eqref{eq:35},
  then the knowledge of $\abs{u(t, x)}$
  uniquely determines 
  the value of the coefficients $A$, $b$ and~$\Re(C)$.
  Indeed, by lemma~\ref{lem:6}
  the imaginary parts $\Im(A)$ and $\Im(b)$
  are determined by the fact
  that $\abs{u(t, x)}$ has a unique maximum 
  at the point \mbox{$t = -\Im(A)$}, \mbox{$x = -\Im(b)$},
  while the real parts $\Re(A) < 0$, $\Re(b)$ and $\Re(C)$
  are determined by the coefficients of the polynomial
  \begin{multline*}
    \Abs{u\tonde{t - \Im(A), x - \Im(b)}}^{-2} = \\
    = 4 \pi^4 \Eu^{-2 \Re(C)} \tonde{%
      \Tonde{(\Re(A))^2 - \abs{\Re(b)}^2 + \abs{x}^2 - t^2}^2 +
      4 \Tonde{\Re(A) t - \Re(b) \cdot x}^2}.
  \end{multline*}
\end{remark}

\medskip
We can repeat the above procedure for the term $u_-^2$
(the only difference is that~$\tau$ must be replaced by~$-\tau$):
\begin{equation} \label{eq:38}
  \Norm{u_-}_{L^4(\R^4)} \le (2 \pi)^{-1/4} \norm{f_-}_{L^2(\R^3)},
\end{equation}
with equality if and only if $f_-$ is of the form~\eqref{eq:35}.

\medskip
For the term $u_+ u_-$,
we observe that by H\"older's inequality we have
\begin{equation} \label{eq:39}
  \Norm{u_+ u_-}_{L^2(\R^4)} \le
  \Norm{u_+}_{L^4(\R^4)} \Norm{u_-}_{L^4(\R^4)} \le
  (2 \pi)^{-1/2} \norm{f_+}_{L^2(\R^3)} \norm{f_-}_{L^2(\R^3)}.
\end{equation}
The first inequality in~\eqref{eq:39} is an equality
if there is a constant $\mu \in \R$
such that $\abs{u_+(t, x)} = \mu \abs{u_-(t, x)}$
for (almost all) $(t, x) \in \R \times \R^3$.
The second inequality in~\eqref{eq:39} is an equality
if $f_+$ and $f_-$ are functions of the form~\eqref{eq:35}.

\medskip
Combining the $L^2$ orthogonality
of the terms $u_+^2$, $u_-^2$ and $u_+ u_-$
(due to the disjointness of the supports
of their Fourier transforms)
with~\eqref{eq:34},~\eqref{eq:38} and~\eqref{eq:39},
we obtain
\begin{multline} \label{eq:40}
  \norm{u}_{L^4}^4 = \norm{(u_+ + u_-)^2}_{L^2}^2 = 
  \norm{u_+}_{L^4}^4 + \norm{u_-}_{L^4}^4
  + 4\norm{u_+ u_-}_{L^2}^2 \le \\
  \le \1{2 \pi} \tonde{\norm{f_+}_{L^2}^4 + \norm{f_-}_{L^2}^4 +
    4 \norm{f_+}_{L^2(\R^3)}^2 \norm{f_-}_{L^2(\R^3)}^2} \le \\
  \le \frac{3}{4 \pi}
  \tonde{\norm{f_+}_{L^2}^2 + \norm{f_-}_{L^2}^2}^2 =
  \frac{3}{16 \pi}
  \norm{(f, g)}_{\dot{H}^{\12} \times \dot{H}^{-\12}}^4,
\end{multline}
where we have used the sharp inequality
\begin{equation*}
  X^2 + Y^2 + 4 X Y \le \frac32 (X + Y)^2, \qquad X,Y \ge 0,
\end{equation*}
for which equality holds if and only if $X = Y$.
This proves that for $n = 3$
the best constant $W(3)$ in~\eqref{eq:3}
is no larger that $(3/(16\pi))^{1/4}$.
The next proposition tells us that maximizers exist
and that the inequalities in~\eqref{eq:40} are sharp;
hence, $W(3) = (3/(16\pi))^{1/4}$.

\begin{proposition} \label{pro:6}
  We have $\norm{u}_{L^4} = (3/(16\pi))^{1/4}
  \norm{(f, g)}_{\dot{H}^{\12} \times \dot{H}^{-\12}}$
  if and only if
  \begin{equation} \label{eq:41}
    \fph(\xi) = \abs{\xi}^{-\12}
    \exp\TOnde{A \abs{\xi} + b \cdot \xi + C}, \qquad
    \fmh(\xi) = \abs{\xi}^{-\12}
    \exp\TOnde{\ovl{A} \abs{\xi} - \ovl{b} \cdot \xi + D},
  \end{equation}
  where $A, C, D \in \C$ and $b \in \C^3$
  with $\abs{\Re(b)} < -\Re(A)$ and $\Re(D) = \Re(C)$.
\end{proposition}

\begin{proof}
  By the above discussion, 
  we have equalities in~\eqref{eq:40} if and only if
  $f_+$, $f_-$ are both functions of the form~\eqref{eq:35}
  and $\abs{u_+(t, x)} = \abs{u_-(t, x)}$
  for all $(t, x) \in \R \times \R^3$.
  Observe that
  if $u_-$ is a $(-)$-wave with data $f_-$, where
  \begin{equation*}
    \fmh(\xi) = \abs{\xi}^{-\12}
    \exp\TOnde{A_- \abs{\xi} + b_- \cdot \xi + C_-},
  \end{equation*}
  then its complex conjugate $\ovl{u_-}$ is a $(+)$-wave with data
  \begin{equation*}
    \ovl{\fmh(-\xi)} = \abs{\xi}^{-\12} \exp\TOnde{%
      \ovl{A_-} \abs{\xi} - \ovl{b_-} \cdot \xi + \ovl{C_-}}.
  \end{equation*}
  By remark~\ref{rem:2}, 
  if two $(+)$-waves with initial data of the form~\eqref{eq:35}
  have the same absolute value at every point of the space-time
  then they must have the same coefficients $A$, $b$ and $\Re(C)$.
\end{proof}

A particular case of~\eqref{eq:41},
corresponding to $A = -1$, $b = 0$, $C = D = \log(2\pi^2)$,
is given by the initial data
\begin{equation} \label{eq:42}
   f_*(x) = \1{1 + \abs{x}^2}, \qquad
   g_*(x) = 0, \qquad x \in \R^3.
\end{equation}

The class of initial data of the form~\eqref{eq:41}
is invariant under the action of the group~$\mL$
described in definition~\ref{def:2}.
The coefficients change according to the following rules:
  \begin{itemize}
  \item space-time translations:
    $(A, b, C, D) \changeto (A + \iC t_0, b + \iC x_0, C, D)$;
  \item isotropic dilations:
    $(A, b, C, D) \changeto
    (A/\lambda, b/\lambda, C - n \log\lambda, D - n \log\lambda)$;
  \item change of scale:
    $(A, b, C, D) \changeto (A, b, C + \log\mu, D + \log\mu)$;
  \item space rotations:
    $(A, b, C, D) \changeto (A, R b, C, D)$;
  \item phase shifts:
    $(A, b, C, D) \changeto (A, b, C + \iC\theta_+, D + \iC\theta_-)$;
  \item Lorentzian boosts:
    \begin{equation*}
      \Tonde{A, (b_1, b'), C, D} \changeto
      \TOnde{A \cosh(a) - b_1 \sinh(a),
        \Tonde{-A \sinh(a) + b_1 \cosh(a), b'}, C, D}.
    \end{equation*}
  \end{itemize}
Hence, after a translation and a phase shift
we can make all coefficients real;
by a rotation we can make $b' = 0$
and by a Lorentzian boost we can make $b_1 = 0$;
then, by a isotropic dilation we can have $A = -1$;
finally a change of scale gives $C = D = \log(2\pi^2)$.
This would correspond to the functions
\begin{equation*}
  \fph(\xi) = \fmh(\xi) =
  -2 \pi^2 \abs{\xi}^{-\12} \exp\Tonde{-\abs{\xi}},
\end{equation*}
which are the Fourier transforms
of the $(+)$ and $(-)$ parts of the initial data~\eqref{eq:42}.
Thus, we have proved that
any maximizer is connected to $(f_*, g_*)$
by the action of~$\mL$.

\section{%
  Wave equation in dimension $n = 2$.}
\label{sec:case-w-n-2}

Consider the case $n = 2$, $p = 6$ for estimate~\eqref{eq:3}.
We decompose $u$ into its $(+)$ and $(-)$ parts
and treat the $L^6$ norm of $u$ as an $L^2$ norm of $u^3$.
By expanding the the products we find
\begin{align*}
  \norm{u}_{L^6}^6 =& \norm{(u_+ + u_-)^3}_{L^2}^2 =
  \norm{u_+^3 + 3 u_+^2 u_- + 3 u_+ u_-^2 + u_-^3}^2 = \\
  =& \norm{u_+^3}^2
  + \norm{u_-^3}^2
  + 9 \norm{u_+^2 u_-}^2
  + 9 \norm{u_+ u_-^2}^2 + \\
  &+ 6 \Re\ip{u_+^3, u_+^2 u_-}
  + 6 \Re\ip{u_+ u_-^2, u_-^3}
  + 18 \Re\ip{u_+^2 u_-, u_+ u_-^2} + \\
  &+ 6 \Re\ip{u_+^3, u_+ u_-^2}
  + 2 \Re\ip{u_+^3, u_-^3}
  + 6 \Re\ip{u_+^2 u_-, u_-^3},
\end{align*}
where here $\norm{\cdot}$ and $\ip{\cdot, \cdot}$ now stand for
the standard norm and inner product in~$L^2(\R \times \R^2)$.
We shall study one term at a time,
but first we compute some integrals which will be needed later.

\begin{lemma} \label{lem:7}
  For $(\tau, \xi) \in \R \times \R^2$ with $\tau > \abs{\xi}$, we define
  \begin{align*}
    I_2(\tau, \xi) &= \int_{(\R^2)^2}
    \ddirac{\tau - \abs{\eta_1} - \abs{\eta_2} \\
      \xi - \eta_1 - \eta_2}
    \frac{\d\eta_1 \d\eta_2}{\abs{\eta_1} \abs{\eta_2}}, \\
    I_3(\tau, \xi) &= \int_{(\R^2)^3}
    \ddirac{\tau - \abs{\eta_1} - \abs{\eta_2} - \abs{\eta_3} \\
      \xi - \eta_1 - \eta_2 - \eta_3}
    \frac{\d\eta_1 \d\eta_2 \d\eta_3}
    {\abs{\eta_1} \abs{\eta_2} \abs{\eta_3}}.
  \end{align*}
  Then we have $I_2(\tau, \xi) = 2 \pi / {\sqrt{\tau^2 - \abs{\xi}^2}}$
  and $I_3(\tau, \xi) = 4 \pi^2$.
\end{lemma}

\begin{proof}
  The fact that $I_2$ and $I_3$ are well defined
  locally integrable functions on $\R \times \R^2$
  follows from considerations similar to the ones made
  at the end of remark~\ref{rem:7}.
  Let us define $\mu$ to be the measure
  $\mu(\tau, \xi) = \abs{\xi}^{-1} \ddirac{\tau - \abs{\xi}}$;
  we have \mbox{$I_2 = \mu \conv \mu$} and $I_3 = \mu \conv \mu \conv \mu$.
  The measure $\mu$ is invariant under proper Lorentz transformations,
  and given $(\tau, \xi)$ such that $\tau > \abs{\xi}$,
  there exists always a proper Lorentz transformation which
  takes $(\tau, \xi)$ to the point $(\tau^*, 0)$
  where $\tau^* = \sqrt{\tau^2 - \abs{\xi}^2}$.
  By lemma~\ref{lem:21},
  it follows that $I_k(\tau, \xi) = I_k(\tau^*, 0)$, for $k = 2, 3$.
  The integral $I_2$ is homogeneous of degree $-1$
  while $I_3$ is homogeneous of degree $0$.
  Hence, for $\tau > \abs{\xi}$ we have
  \begin{equation*}
    I_2(\tau, \xi) = \frac{I_2(1, 0)}{\tau^*} = 
    \1{\tau^*}
    \int_{\R^2} \frac{\ddirac{1 - 2 \abs{\eta}}}{\abs{\eta}^2} \d\eta =
    \frac{2 \pi}{\tau^*} \int_0^\infty \frac{\ddirac{1 - 2 r}}{r} \d r =
    \frac{2 \pi}{\sqrt{\tau^2 - \abs{\xi}^2}}
  \end{equation*}
and
  \begin{multline*}
    I_3(\tau, \xi) = I_3(1, 0) =
    \int_{2 \abs{\eta} \le 1} I_2(1 - \abs{\eta}, -\eta) \frac{\d\eta}{\abs{\eta}} =
    \\ = \int_{2 \abs{\eta} < 1} \frac{2 \pi \d\eta}
    {\Tonde{(1 - \abs{\eta})^2 - \abs{\eta}^2}^{\12} \abs{\eta}}
    = 4 \pi^2 \int_0^{1/2} \frac{\d r}{\sqrt{1 - 2 r}} = 4 \pi^2.
  \end{multline*}
\end{proof}

Let us now begin the proof of the estimate
for the term $\norm{u_+^3}$.
The Fourier transform of $u_+^3$ is
\begin{equation} \label{eq:43}
  \widetilde{u_+^3}(\tau, \xi) =
  \1{(2 \pi)^3} \int_{(\R^2)^3}
  \frac{\fph(\eta_1) \fph(\eta_2) \fph(\eta_3)}
  {\abs{\eta_1}^{\12} \abs{\eta_2}^{\12} \abs{\eta_3}^{\12}}
  \ddirac{%
    \tau - \abs{\eta_1} - \abs{\eta_2} - \abs{\eta_3} \\
    \xi - \eta_1 - \eta_2 - \eta_3}
  \d\eta_1 \d\eta_2 \d\eta_3.
\end{equation}
The support of $\widetilde{u_+^3}$ is contained in the closure of the region
$\mC_{+++} = \graffe{(\tau, \xi): \tau > \abs{\xi}}$.
For each choice of $(\tau, \xi) \in \mC_{+++}$ ,
we denote by $\ip{\cdot, \cdot}_{(\tau, \xi)}$
the $L^2$ inner product associated with the measure
\begin{equation} \label{eq:44}
  \ddirac{\tau - \abs{\eta_1} - \abs{\eta_2} - \abs{\eta_3} \\
    \xi - \eta_1 - \eta_2 - \eta_3}
  \d\eta_1 \d\eta_2 \d\eta_3.
\end{equation}
We can then write
\begin{equation} \label{eq:45}
  \widetilde{u_+^3}(\tau, \xi) = \1{(2 \pi)^3}
  \ip{\fph \otimes \fph \otimes \fph, \abs{\cdot}^{-\12} \otimes
    \abs{\cdot}^{-\12} \otimes \abs{\cdot}^{-\12}}_{(\tau, \xi)},
\end{equation}

\begin{lemma} \label{lem:4}
  For each $(\tau, \xi) \in \mC_{+++}$ we have 
  $\norm{\abs{\cdot}^{-\12} \otimes \abs{\cdot}^{-\12}
    \otimes \abs{\cdot}^{-\12}}_{(\tau, \xi)} = 2 \pi$.
\end{lemma}

\begin{proof}
  The square of the norm we want to compute is the integral $I_3$
  of lemma~\ref{lem:7},
  \begin{equation*}
    \norm{\abs{\cdot}^{-\12} \otimes \abs{\cdot}^{-\12}
      \otimes \abs{\cdot}^{-\12}}_{(\tau, \xi)}^2 = 
    I_3(\tau, \xi) = 4 \pi^2.
  \end{equation*}
\end{proof}

Cauchy-Schwarz's inequality applied to~\eqref{eq:45}
and lemma~\ref{lem:4} give
\begin{multline} \label{eq:46}
  \Norm{u_+^3}^2 = \1{(2 \pi)^3} \Norm{\widetilde{u_+^3}}^2 \le
  \1{(2 \pi)^7} \int_{\R \times \R^2}
  \norm{\fph \otimes \fph \otimes \fph}_{(\tau, \xi)}^2
  \d\tau \d \xi = \\
  = \1{(2 \pi)^7} \norm{%
    \fph \otimes \fph \otimes \fph}_{L^2\tonde{(\R^2)^3}}^2
  = \1{(2 \pi)^7} \Norm{\fph}^6
  = \1{2 \pi} \norm{f_+}^6.
\end{multline}
This time,
equality holds if there exist a function $F: \mC_{+++} \to \C$
such that 
\begin{equation*}
  \Tonde{\fph \otimes \fph \otimes \fph}(\eta_1, \eta_2, \eta_3) =
  F(\tau, \xi)
  \abs{\eta_1}^{-\12} \abs{\eta_2}^{-\12} \abs{\eta_3}^{-\12},
\end{equation*}
for all $(\eta_1, \eta_2, \eta_3)$
in the support of the measure~\eqref{eq:44}.
This means that
\begin{equation*}
  \abs{\eta_1}^{\12} \fph(\eta_1)
  \abs{\eta_2}^{\12} \fph(\eta_2)
  \abs{\eta_3}^{\12} \fph(\eta_3)
  = F\Tonde{\abs{\eta_1} + \abs{\eta_2} + \abs{\eta_3},
    \eta_1 + \eta_2 + \eta_3},
\end{equation*}
for $\eta_1, \eta_2, \eta_3 \in \R^2$.
Examples of such functions are again
$\fph(\xi) = \abs{\xi}^{-\12} \Eu^{-\abs{\xi}}$,
$F(\tau, \xi) = \Eu^{-\tau}$.
More generally, 
by proposition~\ref{pro:5}
all maximizers for the estimate~\eqref{eq:46}
are given by the family
\begin{equation} \label{eq:47}
  \fh(\xi) =
  \abs{\xi}^{-\12} \exp\Tonde{A \abs{\xi} + b \cdot \xi + C},
\end{equation}
with $A, C \in \C$, $b \in \C^2$
and $\abs{\Re(b)} < -\Re(A)$ (in order to have an $L^2$ function).

\begin{lemma} \label{lem:11}
  Let $u$ be the $(+)$-wave corresponding
  to an $L^2$ function of the form~\eqref{eq:47},
  \begin{equation*}
    u(t, x) = \1{(2\pi)^2} \int_{\R^2}
    \exp\TOnde{(A + \iC t) \abs{\xi} + (b + \iC x) \cdot \xi + C}
    \frac{\d\xi}{\abs{\xi}}.
  \end{equation*}
  Then we have the explicit formula
  \begin{multline} \label{eq:48}
    2 \pi \Eu^{-\iC \Im(C)} u\Tonde{t - \Im(A), x - \Im(b)} = \\
    = \frac{\Eu^{\Re(C)}}
    {\sqrt{(\Re(A))^2 - \abs{\Re(b)}^2 + \abs{x}^2 - t^2 +
        2 \iC \tonde{\Re(A) t - \Re(b) \cdot x}}}.
  \end{multline}
\end{lemma}

\begin{proof}
  The integral
  \begin{equation*}
    F(t, x) = \int_{\R^2} \exp\Tonde{t \abs{\xi} + x \cdot \xi}
    \frac{\d\xi}{\abs{\xi}}
  \end{equation*}
  is well defined for $t \in \C$ and $x \in \C^2$
  when $\Re(t) < -\abs{\Re(x)}$.
  For $t \in \R$ and $x \in \R^2$ with $t < -\abs{x}$,
  using polar coordinate we find
  \begin{equation*}
    F(t, x) = \int_0^{2 \pi} \int_0^\infty
    \exp\Tonde{(t + \abs{x} \cos\theta) r} \d r \d\theta =
    \int_0^{2 \pi} \frac{\d\theta}{-t - \abs{x} \cos\theta}
    = \frac{2 \pi}{\sqrt{t^2-\abs{x}^2}}.
  \end{equation*}
  By analytic continuation this formula remains valid
  for complex $t$ and $x$ with $\Re(t) < -\abs{\Re(x)}$.
  Formula~\eqref{eq:48} follows from the identity  
  \begin{equation*}
    \Eu^{-\iC \Im(C)} u\Tonde{t - \Im(A), x - \Im(b)}
    = \frac{\Eu^{\Re(C)}}{(2 \pi)^2} 
    F\Tonde{\Re(A) + \iC t, \Re(b) + \iC x}.
  \end{equation*}
\end{proof}

\begin{remark} \label{rem:3}
  If $u$ is the $(+)$-wave corresponding
  to an $L^2$ function of the form~\eqref{eq:47},
  then the knowledge of $\abs{u(t, x)}$
  uniquely determines 
  the value of the coefficients~$A$,~$b$ and~$\Re(C)$.
  The proof of this fact is similar
  to the one outlined in remark~\ref{rem:2}.
\end{remark}

\medskip
Similarly, for the term $\norm{u_-^3}$ we have
\begin{equation*}
  \Norm{u_-^3}^2 \le \1{2 \pi} \norm{f_-}^6,
\end{equation*}
with equality when $f_-$ takes the form~\eqref{eq:47}.

\medskip
For the term $\norm{u_+^2 u_-}$,
we observe that by H\"older's inequality we have
\begin{equation} \label{eq:49}
  \Norm{u_+^2 u_-}^2 \le
  \Norm{u_+^3}^{\frac43} \Norm{u_-^3}^{\frac23} \le
  \1{2 \pi} \norm{f_+}^4 \norm{f_-}^2.
\end{equation}
The second inequality in~\eqref{eq:49} is an equality
if $f_+$ and $f_-$ are functions of the form~\eqref{eq:47}.
The first inequality in~\eqref{eq:49} is an equality
if there is a constant $\mu \ge 0$
such that $\abs{u_+(t, x)} = \mu \abs{u_-(t, x)}$
for (almost all) $(t, x) \in \R \times \R^2$.

\begin{lemma}
  Let $u_+$ be a $(+)$-wave and $u_-$ be a $(-)$-wave
  corresponding to initial data $f_+$ and $f_-$
  of the form~\eqref{eq:47}.
  If there exists $\mu \ge 0$ 
  such that $\abs{u_+(t, x)} = \mu \abs{u_-(t, x)}$
  for all $t$ and $x$,
  then
  \begin{equation} \label{eq:50}
    \fph(\xi) = \abs{\xi}^{-\12}
    \exp\Tonde{A \abs{\xi} + b \cdot \xi + C}, \qquad
    \fmh(\xi) = \abs{\xi}^{-\12}
    \exp\Tonde{\ovl{A} \abs{\xi} - \ovl{b} \cdot \xi + D},
  \end{equation}
  for some $A, C, D \in \C$ and $b \in \C^2$.
\end{lemma}
The proof of this lemma follows from the same argument
used in the proof of proposition~\ref{pro:6}.

\medskip
Similarly, for the term $\norm{u_+ u_-^2}$ we have
\begin{equation*}
  \Norm{u_+ u_-^2}^2 \le \1{2 \pi} \norm{f_+}^2 \norm{f_-}^4,
\end{equation*}
with equality if and only if
$f_+$ and $f_-$ are functions of the form~\eqref{eq:50}.

\medskip
Let us consider now the term $\Re\ip{u_+^3, u_+^2 u_-}$.
We have
\begin{equation*}
  \Re\ip{u_+^3, u_+^2 u_-} \le
  \abs{\ip{u_+^3, u_+^2 u_-}} \le
  \norm{u_+^3} \norm{u_+^2 u_-} \le
  \1{2 \pi} \norm{f_+}^5 \norm{f_-}.
\end{equation*}
Equality in the second and third inequalities here implies that
$f_+$ and $f_-$ are of the form~\eqref{eq:50},
while we must have $\Im(C) = \Im(D)$
to have equality in the first inequality.

\medskip
Similarly for the terms $\Re\ip{u_+ u_-^2, u_-^3}$
and $\Re\ip{u_+^2 u_-, u_+ u_-^2}$
we have
\begin{equation*}
  \abs{\Re\ip{u_+ u_-^2, u_-^3}} \le
  \1{2 \pi} \norm{f_+} \norm{f_-}^5, \qquad
  \abs{\Re\ip{u_+^2 u_-, u_+ u_-^2}} \le
  \1{2 \pi} \norm{f_+}^3 \norm{f_-}^3,
\end{equation*}
with equality when $f_+$ and $f_-$ are of the form~\eqref{eq:50}
with $\Im(C) = \Im(D)$.

\medskip
The terms $\Re\ip{u_+^3, u_+ u_-^2}$, $\Re\ip{u_+^3, u_-^3}$,
and $\Re\ip{u_+^2 u_-, u_-^3}$ are always zero.
Indeed, the Fourier transform of the cubic terms
$u_+^3$, $u_+^2 u_-$, $u_+ u_-^2$, $u_-^3$
are $L^2$ functions supported on the closures of the regions
\begin{align*}
  \mC_{+++} &=
  \graffe{(\tau, \xi) \in \R \times \R^2: \tau > \abs{\xi}}, \\
  \mC_{++-} &=
  \graffe{(\tau, \xi) \in \R \times \R^2: \tau > -\abs{\xi}}, \\
  \mC_{+--} &=
  \graffe{(\tau, \xi) \in \R \times \R^2: \tau > \abs{\xi}}, \\
  \mC_{---} &= 
  \graffe{(\tau, \xi) \in \R \times \R^2: \tau > -\abs{\xi}},
\end{align*}
respectively,
and the intersections
$\ovl{\mC_{+++}} \intersection \ovl{\mC_{+--}}$,
$\ovl{\mC_{+++}} \intersection \ovl{\mC_{---}}$,
$\ovl{\mC_{++-}} \intersection \ovl{\mC_{---}}$
are sets of measure zero.

\medskip
We put together all the estimates for each single term
and obtain
\begin{multline} \label{eq:51}
  2 \pi \norm{u}_{L^6}^6 \le
  \norm{f_+}^6 + \norm{f_-}^6
  + 9 \norm{f_+}^4 \norm{f_-}^2
  + 9 \norm{f_+}^2 \norm{f_-}^4 + \\
  + 6 \norm{f_+}^5 \norm{f_-}
  + 6 \norm{f_+} \norm{f_-}^5
  + 18 \norm{f_+}^3 \norm{f_-}^3.
\end{multline}

\begin{lemma} \label{lem:12}
  For $X \ge 0$ and $Y \ge 0$ we have
  the sharp polynomial inequality
  \begin{equation*}
    X^6 + Y^6 + 9 X^4 Y^2 + 9 X^2 Y^4 +
    6 X^5 Y + 6 X Y^5 + 18 X^3 Y^3 \le
    \frac{25}{4} \tonde{X^2 + Y^2}^3
  \end{equation*}
  with equality if and only if $X = Y$.
\end{lemma}

\begin{proof}
  By homogeneity we can assume that $Y = 1$.
  Let
  \begin{equation*}
    P(X) = X^6 + 1 + 9 X^4 + 9 X^2 + 6 X^5 + 6 X + 18 X^3, \qquad
    Q(X) = X^2 + 1.
  \end{equation*}
  We want to prove that $4 P(X) \le 25 Q(X)^3$ for $X \ge 0$,
  with equality if and only if $X = 1$.
  Since we have the identity
  \begin{equation*}
    4 P(X) = 4 Q(X)^3 + 24 \tonde{X Q(X)^2 + X^2 Q(X) + X^3},
  \end{equation*}
  our inequality is equivalent to
  \begin{equation*}
    24 \tonde{X Q(X)^2 + X^2 Q(X) + X^3} \le 21 Q(X)^3,
  \end{equation*}
  which reduces to
  \begin{equation*}
    8 \tonde{Z + Z^2 + Z^3} \le 7, \qquad 
    Z = \frac{X}{Q(X)} = \frac{X}{X^2+1} \in [0, 1/2].
  \end{equation*}
  On the interval $[0,1/2]$
  the polynomial $Z + Z^2 + Z^3$ is strictly increasing
  and takes its maximum value $7/8$ when $Z = 1/2$,
  which corresponds to $X = 1$.
\end{proof}

We apply lemma~\ref{lem:12} to~\eqref{eq:51} and finally obtain
\begin{equation} \label{eq:52}
  \norm{u}_{L^6} \le
  \tonde{\frac{25}{8 \pi}}^{\16}
  \tonde{\norm{f_+}_{L^2}^2 + \norm{f_-}^2}^{\12} = 
  \tonde{\frac{25}{64 \pi}}^{\16}
  \norm{(f, g)}_{\dot{H}^{\12}(\R^n) \times \dot{H}^{-\12}(\R^n)},
\end{equation}
which proves that $W(2) \ge (25 / (64 \pi))^{1/6}$.
The next proposition tells us that maximizers exist
and that all the inequalities are sharp;
hence, $W(2) = (25 / (64 \pi))^{1/6}$.

\begin{proposition}
  We have $\norm{u}_{L^6} = (25 / (64 \pi))^{1/6}
  \norm{(f, g)}_{\dot{H}^{\12} \times \dot{H}^{-\12}}$
  if and only if
  \begin{equation} \label{eq:53}
    \fph(\xi) = \abs{\xi}^{-\12}
    \exp\TOnde{A \abs{\xi} + b \cdot \xi + C}, \qquad
    \fmh(\xi) = \abs{\xi}^{-\12}
    \exp\TOnde{\ovl{A} \abs{\xi} - \ovl{b} \cdot \xi + C},
  \end{equation}
  where $A, C \in \C$ and $b \in \C^2$
  with $\abs{\Re(b)} < -\Re(A)$.  
\end{proposition}

\begin{proof}
  {}From the above discussion
  we have equalities in the estimates for each single cubic term
  if and only if $f_+$ and $f_-$ are of the form~\eqref{eq:50}
  with $\Im(C) = \Im(D)$.
  To have equality in~\eqref{eq:52} we also need 
  $\norm{f_+} = \norm{f_-}$
  which in this case implies $\Re(C) = \Re(D)$.
\end{proof}

A particular case of~\eqref{eq:53},
corresponding to $A = -1$, $b = 0$, $C = \log(2\pi)$,
is given by the initial data
\begin{equation} \label{eq:54}
   f_*(x) = \1{\sqrt{1 + \abs{x}^2}}, \qquad
   g_*(x) = 0, \qquad x \in \R^3.
\end{equation}

As we have seen at the end of section~\ref{sec:case-w-n-3},
the class of initial data of the form~\eqref{eq:53}
is invariant under the action of the group $\mL$
and any maximizer is connected to
the functions $(f_*, g_*)$
by the action of $\mL$.

\section{Functional equations}
\label{sec:functional-equations}

In this section we study the functional equations which
characterize the families of maximizers that
we have found in the previous sections.
They are:
\begin{align}
  \label{eq:55}
  & f(x) f(y) f(z) = F\Tonde{x^2 + y^2 + z^2, x + y + z} &
  & \text{for a.e. $(x, y, z) \in \R^3$}, \\
  \label{eq:56}
  & f(x) f(y) = F\Tonde{\abs{x}^2 + \abs{y}^2, x + y} &
  & \text{for a.e. $(x, y) \in \R^2 \times \R^2$}, \\
  \label{eq:57}
  & f(x) f(y) f(z) =
  F\Tonde{\abs{x} + \abs{y} + \abs{z}, x + y + z}, &
  & \text{for a.e. $(x, y, z) \in (\R^2)^3$}, \\
  \label{eq:58}
  & f(x) f(y) = F\Tonde{\abs{x} + \abs{y}, x + y}, &
  & \text{for a.e. $(x, y) \in \R^3 \times \R^3$},
\end{align}
where $f$ and $F$
are unknown complex valued measurable functions
and the identities are supposed to hold almost everywhere
with respect to the Lebesgue measure.
We are going to show that locally integrable solutions
to these equations are actually smooth functions;
this is a general principle which holds for a large class
of functional equation
(actually even assuming only measurability implies continuity,
see the work of A.\ J\'arai in \cite{Jar2003}),
but for the sake of completeness we include a direct proof
adapted to our equations.
Once the smoothness of $f$ and $F$ is established,
it is not difficult to solve the equation using 
geometric or differential methods.
It turns out that in all cases the function $F$
must be an exponential function of the form
\begin{equation*}
  F(t, x) = \exp\Tonde{A t + b \cdot x + C}.
\end{equation*}

A simpler model for the above functional equations
is provided by the exponential law, $f(x) f(y) = f(x + y)$,
which is one of the four basic Cauchy functional equations.
We refer the reader to~\cite{Acz1966} for a general introduction
to the subject of functional equations.
Here, we only require the following result
which is a simple exercise in real analysis.

\begin{lemma} \label{lem:10}
  Let $\Omega$ be an open subset of $\R^n$ 
  such that $x + y \in \Omega$ whenever $x, y \in \Omega$.
  Let $f: \Omega \to \C$ be a non-trivial locally integrable
  solution of the Cauchy functional equation
  \begin{equation} \label{eq:59}
    f(x) f(y) = f(x + y), \qquad
    \text{for a.e.\ $(x, y) \in \Omega \times \Omega$}.
  \end{equation}
  Then there exists a vector $b \in \C^n$ such that 
  $f(x) = \exp(b \cdot x)$ for a.e.\ $x \in \Omega$.
\end{lemma}

\begin{proof}
  Let $Q$ be cube contained in $\Omega$ such that
  $\int_Q f(y) \d y \ne 0$.
  If we integrate~\eqref{eq:59} with respect to $y \in Q$
  we obtain that $f(x)$ must coincide (almost everywhere)
  with the continuous function
  \begin{equation*}
    x \mapsto \frac{\int_Q f(x+y) \d y}{\int_Q f(y) \d y}.
  \end{equation*}
  If $f$ is continuous the above function is differentiable.
  Hence, we may assume that $f$ is differentiable.
  Fix $y_0 \in \Omega$ and let $b = (\nabla f(y_0))/f(y_0)$.
  If we differentiate~\eqref{eq:59} with respect to $y$
  and set $y = y_0$ we obtain the differential equation
  \begin{equation*}
    \nabla f(x + y_0) = f(x) \nabla f(y_0) =
    f(x) f(y_0) b = f(x + y_0) b,
  \end{equation*}
  whose non trivial solutions have the form
  $f(z) = \exp(b \cdot z + C)$ for some constant $C \in \C$.
  Substituting this expression for $f$ into~\eqref{eq:59},
  we obtain that $\exp(C)$ must be $1$.
\end{proof}

As was done in the lemma,
regularity properties of solutions to functional equations
can be obtained by (partial) integration of the equation.
The following lemmata,
although not expressed in their most general form,
are what we need to deduce continuity from local integrability
in our equations.

\begin{lemma} \label{lem:19}
  Let $A, B, \Omega$ be open subsets of $\R^n$.
  Let $f \in L^p_{\rm loc}(\Omega)$
  and let \mbox{$\varphi: A \times B \to \Omega$} be a smooth map such that
  \begin{equation} \label{eq:60}
    \det \abs{\pder{\varphi}{y}(x, y)} \ne 0, \qquad
    \text{for $(x, y) \in A \times B$}.
  \end{equation}
  Let $K$ be a compact subset of $B$.
  For each $x \in A$,
  let $g_x: B \to \Omega$ be the function $g_x(y) = f\tonde{\varphi(x, y)}$.
  Then the map $x \mapsto g_x$
  is a continuous application from $A$ to $L^p(K)$.
\end{lemma}

\begin{proof}
  The case of $f$ continuous is immediate.
  The general case follows by density.
  The condition on the partial Jacobian of $\varphi$
  is enough to apply, at least locally, 
  the change of variable $y \changeto z = \varphi(x, y)$
  in the integration over $K$.
  We leave the details to the reader.
\end{proof}

\begin{remark}
  Lemma~\ref{lem:19} does not hold if we remove condition~\eqref{eq:60}.
  For example, if $A = B = \Omega$ and $\varphi(x, y) = x$
  then we would have $g_x(y) = f(x)$
  and the map $x \mapsto \norm{g_x}_{L^p(K)} = \abs{K}^{1/p} \abs{f(x)}$
  may not be continuous.
\end{remark}

\begin{lemma} \label{lem:15}
  Let $A, B, \Omega$ be open subsets of $\R^n$.
  Let $f_j \in L^{p_j}_{\rm loc}(\Omega)$, $j = 1, 2, \dots, N$,
  with $\sum_j 1/p_j = 1$.
  Let $\varphi_1, \varphi_2, \dots, \varphi_N: A \times B \to \Omega$
  be smooth maps such that
  \begin{equation*}
    \det \abs{\pder{\varphi_j}{y}(x, y)} \ne 0, \qquad
    (x, y) \in A \times B, \quad
    j = 1, 2, \dots, N.
  \end{equation*}
  Let $\psi: A \times B \to \C$ be a continuous function.
  Let $K$ be a compact subset of $B$.
  Then the function
  \begin{equation*}
    F(x) = \int_K \prod_{j = 1}^N f_j\tonde{\varphi_j(x, y)} \psi(x,y) \d y,
    \qquad x \in A,
  \end{equation*}
  is continuous.
\end{lemma}

\begin{proof}
  The continuity of $F$ follows from the previous lemma
  and the continuity of the functional
  \begin{equation*}
    \tonde{x, g_1, g_2, \dots g_N} \mapsto
    \int_K \prod_j g_j(y) \psi(x,y) \d y
  \end{equation*}
  defined on $A \times \prod_j L^{p_j}_{\rm loc}(\Omega)$.
\end{proof}

\begin{proposition} \label{pro:2}
  Let $\Omega$ be an open subset of $\R^n \times \R^n$
  such that the section
  \mbox{$\Omega_x = \graffe{y \in \R^n: (x, y) \in \Omega}$}
  is dense in $\R^n$ for each $x \in \R^n$.
  Let $P, Q: \Omega \to \R^n$ be smooth maps such that
  \begin{equation} \label{eq:61}
    \det \abs{\pder{P}{y}(x, y)} \ne 0, \qquad
    \det \abs{\pder{Q}{y}(x, y)} \ne 0, \qquad
    \text{for $(x, y) \in \Omega$}.
  \end{equation}
  If $f: \R^n \to \C$ is a locally integrable solution
  of the functional equation
  \begin{equation} \label{eq:62}
    f(x) f(y) = f\Tonde{P(x, y)} f\Tonde{Q(x, y)}, \quad
    \text{for a.e. $(x, y) \in \Omega$},
  \end{equation}
  then $f$ is continuous.
\end{proposition}

\begin{proof}
  We may assume that $f$ is non trivial.
  Let $g(x) = \sqrt{\abs{f(x)}}$.
  When \mbox{$f \in L^1_{\rm loc}(\R^n)$} we have $g \in L^2_{\rm loc}(\R^n)$.
  Hence, by~\eqref{eq:61}, it follows that for every $x \in \R^n$
  the function $y \mapsto g(P(x, y)) g(Q(x, y))$
  is locally integrable on $\Omega$.
  Fix $x_0 \in \R^n$ and choose a compact domain $D$ in $\Omega_{x_0}$
  such that $\int_D f \ne 0$.

  We integrate the square root of the absolute value
  of equation~\eqref{eq:62}
  with respect to $y$ over the domain $D$ and obtain
  \begin{equation*}
    g(x) \int_D g(y) \d y = \int_D g\Tonde{P(x, y)} g\Tonde{Q(x, y)} \d y.
  \end{equation*}
  By lemma~\ref{lem:15},
  the right hand side is a continuous function of $x$
  for $x$ in a neighborhood of $x_0$.
  Since $\int_D g \ne 0$, we obtain that $g$ is continuous in $x_0$.
  This proves that $\abs{f}$ is a continuous function.
  In particular it follows that $f \in L^2_{\rm loc}$.
  We can bootstrap the argument:
  if we integrate~\eqref{eq:62}
  with respect to $y$ over the domain $D$, we obtain
  \begin{equation*}
    f(x) \int_D f(y) \d y = 
    \int_D f\Tonde{P(x, y)} f\Tonde{Q(x, y)} \d y,
  \end{equation*}
  from which it follows that $f$ is continuous.
\end{proof}

In the following subsections
we will study in detail each of our four functional equations.
In each case we will adopt the following strategy:
\begin{enumerate}
\item Local integrable solutions are continuous.
\item Nontrivial continuous solutions never vanish.
\item Continuous solutions which never vanish are of exponential form.
\end{enumerate}

\begin{remark} \label{rem:6}
  It is interesting to observe that
  there are functional equations 
  which formally look
  very similar to the ones we are considering,
  but
  for which each of the above three steps fail.
  Take for instance the functional equation
  \begin{equation*}
    f(x) f(y) = F\Tonde{x^2 + y^2, x + y} \text{ for a.e. $(x, y) \in \R^2$}.
  \end{equation*}
  In this case,
  given \emph{any} function $f$
  we can always construct a solution by setting
  \begin{equation*}
    F(s, t) =
    f\TOnde{\frac{t + \sqrt{2 s - t^2}}2}
    f\TOnde{\frac{t - \sqrt{2 s - t^2}}2}.
  \end{equation*}
  Even if the function $f$ is locally integrable,
  it does not need to be an exponential,
  it can vanish on any set and does not need to be continuous.
\end{remark}

\subsection{The equation~\eqref{eq:55}.}

As in section~\ref{sec:case-n-1}
we let $\mP_1 = \graffe{(s, t) \in \R \times \R: 3s > t^2}$.

\begin{lemma} \label{lem:9}
  Let $f: \R \to \C$ and $F: \ovl{\mP_1} \to \C$ be functions which 
  solve equation~\eqref{eq:55}.
  If $f$ is locally integrable then $f$ and $F$ are continuous functions.
\end{lemma}

\begin{proof}
  We first prove that
  $f$ locally integrable implies $F$ locally integrable.
  Indeed, 
  \begin{multline*}
    \iiint_{\abs{x}^2 + \abs{y}^2 + \abs{z}^2 \le R^2}
    \Abs{f(x) f(y) f(z)} \d x\d y\d z = 
    \int_{\substack{v \in \R^3 \\ \abs{v} \le R}}
    \Abs{F\Tonde{\abs{v}^2, v \cdot (1, 1, 1)}} \d v = \\
    = 2 \pi \int_0^R \int_{-1}^1
    \abs{F\Tonde{r^2, \sqrt3 r u}} r^2 \d u\d r =
    \frac{\pi}{\sqrt3} \iint_{\frac{t^2}3 \le s \le R^2}
    \abs{F(s,t)} \d s \d t,
  \end{multline*}
  for any $R > 0$.
  Moreover, using the change of variables
  \begin{equation*}
    (y, z) \changeto (s, t) = \Phi(y, z) = (y^2 + z^2, y + z)
  \end{equation*}
  from the region $\graffe{y > z}$ to $\graffe{s > t^2/2}$,
  with $\d s\d t = 2(y-z) \d y\d z$,
  we have
  \begin{multline*}
    f(x) \int_\Omega f(y) f(z) (y-z) \d y \d z
    = \int_\Omega F\tonde{x^2 + y^2 + z^2, x + y + z}
    (y-z) \d y \d z = \\
    = \12 \int_{\Phi(\Omega)} F\tonde{x^2 + s, x + t} \d s \d t.
  \end{multline*}
  for any bounded domain $\Omega \subseteq \graffe{y > z}$.
  The local integrability of $F$ implies that the function
  $x \mapsto \int_{\Phi(\Omega)} F(x^2 + s, x + t) \d s\d t$
  is continuous.
  We choose $\Omega$ so that
  the integral $\int_\Omega f(y) f(z) (y-z) \d y \d z$ is not zero
  and it follows that $f(x)$ is continuous.

  The continuity of $F$ comes easily from the equation 
  and the continuity of $f$, since we have
  \begin{equation*}
    F(s, t) = f(0)
    f\TOnde{\frac{t + \sqrt{2 s - t^2}}2}
    f\TOnde{\frac{t - \sqrt{2 s - t^2}}2}.
  \end{equation*}
\end{proof}

\begin{remark} \label{rem:1}
  We can write~\eqref{eq:55} as
  \begin{equation*}
    \tonde{f \otimes f \otimes f} (v) = F(s, t), \qquad
    v \in \R^3, \qquad
    s = \abs{v}^2, \qquad
    t = v \cdot (1, 1, 1),
  \end{equation*}
  which shows that the tensor product $f \otimes f \otimes f$,
  as a function on $\R^3$,
  is constant along any circle $\Gamma(s, t)$
  obtained as the intersection
  between the sphere of radius~$\sqrt{s}$ centered at the origin
  and the plane orthogonal to the vector $(1, 1, 1)$
  passing through the point $(t/3, t/3, t/3)$.
\end{remark}

\begin{lemma} \label{lem:8}
  If $f$ and $F$ are continuous functions which 
  solve equation~\eqref{eq:55}
  and $f$ vanishes at one point
  then $f$ and $F$ vanish everywhere.
\end{lemma}

\begin{proof}
  By continuity it is enough to prove that
  the set of all points where $f$ vanishes is open.
  Suppose $f$ vanishes at a point $x_*$.
  If $f$ is not identically zero then
  there exists some open set $A$ in $\R$
  on which $f$ never vanishes.
  Let $y_*$ and $z_*$ be two distinct points in $A$.
  Consider the map
  \begin{equation*}
    (y, z) \mapsto (s, t) = (x_*^2 + y^2 + z^2, x_* + y + z).
  \end{equation*}
  Its Jacobian determinant at the point $(y_*, z_*)$
  is $2 \abs{y_* - z_*} \ne 0$.
  Hence, the map
  is invertible from a neighborhood $U \subset A \times A$ of $(y_*, z_*)$
  to a neighborhood $V$
  of $(s_*, t_*) = (x_*^2 + y_*^2 + z_*^2, x_* + y_* + z_*)$.
  It follows that,
  for each $x$ sufficiently close to $x_*$ so that
  $(x^2 + y_*^2 + z_*^2, x + y_* + z_*)$ lies in $V$,
  there exists a pair $(y, z) \in U$ such that
  \begin{equation*}
    (x^2 + y_*^2 + z_*^2, x + y_* + z_*) = (x_*^2 + y^2 + z^2, x_* + y + z).
  \end{equation*}
  Using the functional equation~\eqref{eq:55} we obtain
  \begin{equation*}
    f(x) f(y_*) f(z_*) = f(x_*) f(y) f(z) = 0,
  \end{equation*}
  and we know that $f(y_*) f(z_*) \ne 0$.
  Hence, $f(x)$ vanishes for $x$ in a neighborhood of~$x_*$.
\end{proof}

\begin{proposition} \label{pro:1}
  If $f: \R \to \C$ and $F: \ovl{\mP_1} \to \C$ 
  are non-trivial locally integrable functions
  which satisfy the functional equation~\eqref{eq:55}
  for all $x, y, z \in \R$,
  then there exists complex constants $A$, $B$, $C$
  such that 
  \begin{equation*}
    f(x) = \exp\Tonde{A x^2 + B x + C}, \qquad
    F(t, x) = \exp\tonde{A t +  B x + 3C},
  \end{equation*}
  for (almost) all $(t, x) \in \mP_1$.
\end{proposition}

\begin{proof}
  By lemma~\ref{lem:9},
  we may assume that $f$ and $F$ are continuous.
  By lemma~\ref{lem:8},
  we may assume that $f$ and $F$ never vanishes.
  We define
  \begin{align*}
    & g(x) = \frac{f(x)}{f(-x)}, &
    & G(s, t) = \frac{F(s, t)}{F(s, -t)}, &
    & h(x) = f(x) f(-x), &
    & H(s, t) = F(s, t) F(s, -t).
  \end{align*}

  The function $g$ corresponds to the \emph{odd} component of $f$,
  $g(x) g(-x) = 1$, $g(0) = 1$,
  and satisfies the same equation
  \begin{equation*}
    g(x) g(y) g(z) = G\tonde{x^2 + y^2 + z^2, x + y + z}, \qquad
    x, y, z \in \R.
  \end{equation*}
  In particular,
  \begin{equation*}
    G(2s, 0) = g\tonde{\sqrt{s}} g\tonde{-\sqrt{s}} g(0) = 1,
  \end{equation*}
  for all $s \ge 0$.
  It follows that
  \begin{equation*}
    g(x) g(y) = g(x + y) g(-x-y) g(x) g(y) = 
    g(x + y) G\tonde{(x+y)^2 + x^2 + y^2, 0} = g(x + y).
  \end{equation*}
  By lemma~\ref{lem:10},
  $g$ must be an exponential function of the form
  $g(x) = \exp(2 B x)$ for some complex constant $B$.

  The function $h$ corresponds to the \emph{even} component of $f$,
  $h(x) = h(-x)$ and satisfies the equation
  \begin{equation*}
    h(x) h(y) h(z) =
    H\tonde{x^2 + y^2 + z^2, x \pm y \pm z}, \qquad
    x, y, z \in \R,
  \end{equation*}
  for any combination of signs.
  By the same argument used in remark~\ref{rem:1},
  we have that $h \otimes h \otimes h$
  is constant along circles obtained by
  intersecting spheres centered at the origin
  with planes perpendicular
  to the four vectors $(1, \pm 1, \pm 1)$.
  It follows that $h \otimes h \otimes h$ must be constant
  on any sphere centered at the origin.
  In fact, any two points at the same distance from the origin
  can be connected by a finite sequence
  of arcs of the above circles.
  This means that there exists
  a function $\varphi: \R_+ \to \C$ such that
  \begin{equation*}
    \frac{h(x) h(y) h(z)}{h(0)^3} =
    \varphi\tonde{x^2 + y^2 + z^2}, \qquad
    x, y, z \in \R.
  \end{equation*}
  In particular, for $s \ge 0$ and $t \ge 0$, we have
  \begin{equation*}
    \varphi(s) \varphi(t) = 
    \frac{h\tonde{\sqrt{s}} h(0)^2}{h(0)^3} \cdot  
    \frac{h\tonde{\sqrt{t}} h(0)^2}{h(0)^3} =
    \frac{h\tonde{\sqrt{s}} h\tonde{\sqrt{t}} h(0)}{h(0)^3} =
    \varphi(s+t).
  \end{equation*}
  By lemma~\ref{lem:10},
  $\varphi$ must be an exponential function of the form
  $\varphi(s) = \exp(2 A s)$ for some complex constant $A$.
  Hence, $h(x) = h(0) \exp\tonde{2 A x^2}$.

  We conclude the proof of the lemma by observing that
  \begin{equation*}
    f(x)^2 = g(x) h(x) = \exp(2 B x) h(0) \exp\tonde{A x^2} =
    \exp\tonde{2 A x^2 + 2 B x + 2 C},
  \end{equation*}
  where $C$ is a complex constant such that $f(0) = \Eu^C$.
\end{proof}

\subsection{The equation~\eqref{eq:56}.}

If we integrate equation~\eqref{eq:56}
with respect to $y$ on a domain of $\R^2$,
we cannot apply directly the regularity results of lemma~\ref{lem:15},
because the domain of $F$ is a region in $\R^3$
and the image of the application
\mbox{$y \mapsto \Tonde{\abs{x}^2 + \abs{y}^2, x + y}$}
is a set of measure zero in $\R^3$.
To overcome this difficulty
we exploit the geometric invariance properties of the equation.

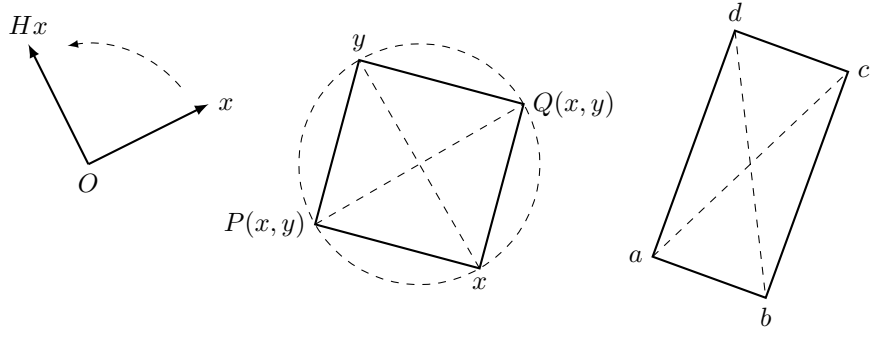
\begin{figure}
  \centering
  \begin{tikzpicture}[scale=0.8, >=latex]
    \begin{scope}
      \draw[->, thick]
      (0, 0) node[below] {$O$} -- (2, 1) node[right] {$x$};
      \draw[->, thick]
      (0, 0) -- (-1, 2) node[above] {$Hx$};
      \draw[->, dashed]
      (0, 0) +(40:2cm) arc (40:100:2cm);
    \end{scope}
    \begin{scope}[xshift=5.5cm, rotate=30]
      \draw[dashed]
      (0, 0) circle (2cm)
      (0, -2) node[below] {$x$} -- (0, 2) node[above] {$y$}
      (-2, 0) node[left] {$P(x,y)$} -- (2, 0) node[right] {$Q(x,y)$};
      \draw[thick]
      (0, -2) -- (2, 0) -- (0, 2) -- (-2, 0) -- cycle;
    \end{scope}
    \begin{scope}[xshift=11cm, rotate=-20]
      \draw[dashed]
      (-1, -2) -- (1, 2)
      (-1, 2) -- (1, -2);
      \draw[thick]
      (-1, -2) node[left] {$a$} -- (1, -2) node[below] {$b$}
      -- (1, 2) node[right] {$c$} -- (-1, 2) node[above] {$d$} -- cycle;
    \end{scope}
  \end{tikzpicture}
  \caption{The functions $H$ (left) and $P, Q$ (center)
    described in remark~\ref{rem:4}
    and a rectangle (right) used in remark~\ref{rem:5}}
  \label{fig:1}
\end{figure}

\begin{remark} \label{rem:4}
  Let $I: \R^2 \to \R^2$ be the identity on $\R^2$ 
  and let $H: \R^2 \to \R^2$, $H(x_1, x_2) = (-x_2, x_1)$,
  be a counterclockwise rotation by $\pi/2$ of the plane.
  Given two points $x$ and $y$ in $\R^2$
  the functions
  \begin{align*}
    P(x, y) &= \frac{x + y}{2} + H\TOnde{\frac{x - y}{2}} =
    \frac{I+H}{2} x + \frac{I-H}{2} y, \\
    Q(x, y) &= \frac{x + y}{2} - H\TOnde{\frac{x - y}{2}} =
    \frac{I-H}{2} x + \frac{I+H}{2} y,
  \end{align*}
  determine two other points $p = P(x, y)$ and $q = Q(x, y)$
  so that $x, y$ and $p, q$ are the opposite vertices of a square
  (see figure~\ref{fig:1}) and we have
  \begin{equation} \label{eq:63}
    p + q = x + y, \qquad
    \abs{p}^2 + \abs{q}^2 = \abs{x}^2 + \abs{y}^2.
  \end{equation}
  Moreover, the linear map $(x, y) \mapsto (p, q)$
  is an isometry on $\R^2 \times \R^2$.
  By property~\eqref{eq:63},
  it follows that the equation~\eqref{eq:56}
  implies the equation
  \begin{equation} \label{eq:64}
    f(x) f(y) = f\Tonde{P(x, y)} f\Tonde{Q(x, y)}, \qquad
    \text{for a.e. $(x, y) \in \R^2 \times \R^2$}.
  \end{equation}
  Observe also that $\de P / \de y = (I - H) / 2$
  and $\de Q / \de y = (I + H) / 2$ are non singular matrices.
\end{remark}

As in section~\ref{sec:case-n-2}
we let $\mP_2 = \graffe{(s, v) \in \R \times \R^2: 2s > \abs{v}^2}$.

\begin{lemma} \label{lem:13}
  Let $f:\R^2 \to \C$ and $F: \ovl{\mP_2} \to \C$
  be solutions of equation~\eqref{eq:56}.
  If $f$ is locally integrable
  then $f$ and $F$ are continuous functions.
\end{lemma}

\begin{proof}
  By remark~\ref{rem:4},
  the lemma becomes a corollary of proposition~\ref{pro:2}.
\end{proof}

\begin{lemma} \label{lem:5}
  If $f$ is a continuous solution of equation~\eqref{eq:64}
  and $f$ vanishes at one point then $f$ vanishes everywhere.
\end{lemma}

\begin{proof}
  Let $y \in \R^2$.
  If $f$ vanishes at the point $x_0$,
  using equation~\eqref{eq:64} we must have that
  $f$ vanishes at the point $x_1$,
  where $x_1$ is either $P(x_0, y)$ or $Q(x_0, y)$,
  and $\abs{x_1 - y} = (1/\sqrt2) \abs{x_0 - y}$.
  By iterating this argument,
  we can construct a sequence of points $x_n$
  such that $f(x_n) = 0$ and $\lim_n x_n = y$.
  By continuity it follows that~$f(y) = 0$.
\end{proof}

\begin{remark} \label{rem:5}
  It follows from equation~\eqref{eq:56}
  that $f$ solves the rectangular functional equation
  \begin{equation*}
    f(a) f(c) = f(b) f(d),
  \end{equation*}
  whenever the points $a, c$ and $b, d$
  are the opposite vertices of a rectangle (see figure~\ref{fig:1}).
  Indeed, when $a - b = d - c$ and $a - b \perp c - b$,
  we have
  \begin{equation*}
    0 = (a - b) \cdot (c - b) =
    a \cdot c - (a + c) \cdot b + \abs{b}^2 = 
    a \cdot c - (b + d) \cdot b + \abs{b}^2 = 
    a \cdot c - d \cdot b.
  \end{equation*}
  Hence, $a \cdot c = b \cdot d$ and
  $\abs{a}^2 + \abs{c}^2 = \abs{a + c}^2 - 2 a \cdot c = 
  \abs{b + d}^2 - 2 b \cdot d = \abs{b}^2 + \abs{d}^2$.
\end{remark}

\begin{proposition} \label{pro:3}
  If $f: \R^2 \to \C$ and $F: \ovl{\mP_2} \to \C$ 
  are nontrivial locally integrable functions
  which satisfy the functional equation~\eqref{eq:56}
  then there exists constants
  \mbox{$A \in \C$}, \mbox{$b \in \C^2$}, $C \in \C$
  such that 
  \begin{equation*}
    f(x) = \exp\Tonde{A \abs{x}^2 + b \cdot x + C}, \qquad
    F(t, x) = \exp\Tonde{A t +  b \cdot x + 2C},
  \end{equation*}
  for (almost) all $(t, x) \in \mP_2$.
\end{proposition}

In the proof of the proposition we follow a geometric construction
which is an adaptation of the one for odd orthogonally additive mappings
found in~\cite{Rat1985}.

\begin{proof}
  By lemma~\ref{lem:13},
  we may assume that $f$ and $F$ are continuous.
  By lemma~\ref{lem:5},
  we may assume that $f$ and $F$ never vanishes.
  We define
  \begin{align*}
    & g(x) = \frac{f(x)}{f(-x)}, &
    & G(t, z) = \frac{F(t, z)}{F(t, -z)}, &
    & h(x) = f(x) f(-x), &
    & H(t, z) = F(t, z) F(t, -z).
  \end{align*}
  The function $g$ corresponds to the \emph{odd} component of $f$,
  $g(x) g(-x) = 1$ and $g(0) = 1$.
  By remark~\ref{rem:5} we know
  it satisfies the rectangular equation $g(a) g(c) = g(b) g(d)$,
  whenever the points $a, c$ and $b, d$
  are the opposite vertices of a rectangle.
  Given two vectors $x$ and $y$ in $\R^2$
  it is always possible to find a third vector $z$ such that
  $z \perp x + y$ and $x + z \perp y - z$.
  Let $p$ and $-p$ be the components of $x$ and $y$ perpendicular to~$x + y$.
  \begin{figure}
    \centering
    \begin{tikzpicture}[scale=0.18, >=latex]
      \begin{scope}
        \draw[dotted]
        (-13, 0) -- (13, 0);
        \draw[->, thick] (0, 0) -- (7, 9);
        \draw[->, thick] (0, 0) -- (-7, 16);
        \draw[->, thick] (0, 0) -- (0, 25);
        \draw[dashed]
        (7, 9) -- (0, 25) -- (-7, 16);
        \draw
        (0, 0) node[below] {$O$} --
        (12, 9) node[right] {$x{+}z$} --
        (12, 0) node[below right] {$p{+}z$} --
        (7, 0) node[below] {$p$} --
        (7, 9) node[left] {$x$} --
        (12, 9) --
        (0, 25) node[above] {$x{+}y$} --
        (-12, 16) node[left] {$y{-}z$} --
        (-12, 0) node[below] {$-p{-}z$} --
        (-7, 0) node[below] {$-p$} --
        (-7, 16) node[right] {$y$} --
        (-12, 16) -- cycle;
      \end{scope}
      \begin{scope}[shift={(30, 2)}]
        \draw (-10,0) node[below] {$-q$} rectangle (0, 15) node[above] {$p$}
        rectangle (10, 0) node[below] {$q$};
        \draw[->, thick] (0, 0) node[below] {$O$}
        -- (-10, 15) node[above] {$-y$};
        \draw[->, thick] (0, 0) -- (10, 15) node[above] {$y$};
      \end{scope}
    \end{tikzpicture}
    \caption{Constructions for the function $g$ (left)
      and the function $h$ (right)
      in the proof of proposition~\ref{pro:3}}
    \label{fig:2}
  \end{figure}
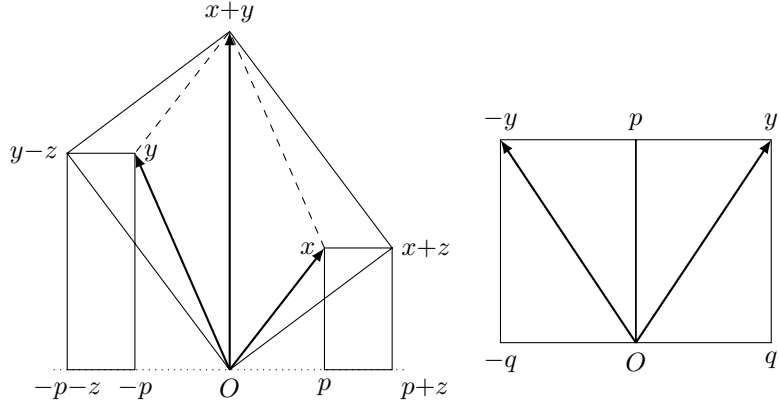
  Consider the three rectangles formed by
  $(0, x+z, x+y, y-z)$, $(x, x+z, p+z, p)$ and $(y, y-z, -p-z, -p)$
  (see figure~\ref{fig:2});
  using the rectangular equation we have
  \begin{align*}
    g(x+y) g(0) &= g(x+z) g(y-z), \\
    g(x) g(p+z) &= g(x+z) g(p), \\
    g(y) g(-p-z) &= g(y-z) g(-p);
  \end{align*}
  and using the parity properties of $g$ we obtain
  \begin{multline*}
    g(x+y) = g(x+y) g(0) = g(x+z) g(y-z) = g(x+z) g(p) g(y-z) g(-p) = \\
    = g(x) g(p+z) g(y) g(-p-z) = g(x) g(y).
  \end{multline*}
  By lemma~\ref{lem:10},
  $g$ must be an exponential function of the form
  $g(x) = \exp(2 b \cdot x)$ for some complex vector $b \in \C^2$.

  The function $h$ corresponds to the \emph{even} component of $f$,
  $h(x) = h(-x)$, and satisfies the rectangular equation
  $h(a) h(c) = h(b) h(d)$,
  whenever the points $a, c$ and $b, d$
  are the opposite vertices of a rectangle.
  Given two points $x$ and $y$ with $\abs{x} = \abs{y}$,
  let $p = (x+y)/2$, $q = (x-y)/2$ and consider the rectangles
  $(0, p, x, q)$, $(0, p, y, -q)$ (see figure~\ref{fig:2});
  by the rectangular equation and the parity of $h$ we have
  \begin{equation*}
    h(x) h(0) = h(p) h(q) = h(p) h(-q) = h(y) h(0).
  \end{equation*}
  Hence, $\abs{x} = \abs{y}$ implies $h(x) = h(y)$.
  This means that $h$ is spherically symmetric
  and there exists a function $\varphi: \R_+ \to \C$ such that
  $h(x) = h(0) \varphi\Tonde{\abs{x}^2}$, for $x \in \R^2$.
  Given $s \ge 0$ and $t \ge 0$,
  let $x$ and $y$ be two points in $\R^2$
  such that $\abs{x}^2 = s$, $\abs{y}^2 = t$ and $x \perp y$;
  by the Pythagorean theorem we have $\abs{x+y}^2 = s + t$.
  It follows that
  \begin{equation*}
    \varphi(s) \varphi(t) = \frac{h(x) h(y)}{h(0)^2} =
    \frac{h(x + y) h(0)}{h(0)^2} = \varphi(s+t).
  \end{equation*}
  By lemma~\ref{lem:10},
  $\varphi$ must be an exponential function of the form
  $\varphi(s) = \exp(2 A s)$ for some complex constant $A$.
  Hence, $h(x) = h(0) \exp\Tonde{2 A \abs{x}^2}$.

  We conclude the proof of the lemma by observing that
  \begin{equation*}
    f(x)^2 = g(x) h(x) = \exp(2 b \cdot x) h(0) \exp\Tonde{A \abs{x}^2} =
    \exp\tonde{2 A x^2 + 2 b \cdot x + 2 C},
  \end{equation*}
  where $C$ is a complex constant such that $f(0) = \Eu^C$.
\end{proof}

\subsection{The equation~\eqref{eq:57}.}

As in section~\ref{sec:case-w-n-2}
we let $\mC_{+++} = \graffe{(t, v) \in \R \times \R^2: t > \abs{v}}$.

\begin{lemma} \label{lem:14}
  Let $f: \R^2 \to \C$ and $F: \ovl{\mC_{+++}} \to \C$
  be functions which solve equation~\eqref{eq:57}.
  If $f$ is locally integrable then $f$ and $F$ are continuous functions.
\end{lemma}

\begin{proof}
  Suppose first that $f \in L^p_{\rm loc}(\R^2)$ for some $p > 2$.
  Using the results of lemma~\ref{lem:7} we can see that 
  $F \in L^1_{\rm loc}(\mC_{+++})$; indeed,
  \begin{multline*}
    \int_{\abs{v}^2 \le t \le R} \abs{F(t, v)} \d v\d t = 
    \1{8 \pi^2} \int_{\abs{v}^2 \le t \le R} \abs{F(t, v)} I_3(t, v) \d v\d t = \\
    = \1{8 \pi^2} \int_{t \le R}
    \frac{\Abs{F\Tonde{\abs{x} + \abs{y} + \abs{z}, x + y + z}}}
    {\abs{x} \abs{y} \abs{z}}
    \ddirac{t - \abs{x} - \abs{y} - \abs{z} \\ v - x - y - z}
    \d x \d y \d z \d v \d t = \\
    = \1{8 \pi^2} \int_{\abs{x} + \abs{y} + \abs{z} \le R}
    \frac{\abs{f(x) f(y) f(z)}}{\abs{x} \abs{y} \abs{z}}
    \d x \d y \d z  \le \1{8 \pi^2}
    \TOnde{\int_{\abs{x} \le R} \frac{\abs{f(x)}}{\abs{x}} \d x}^3 \le \\
    \le C \tonde{R^{1 - 2/p}}^3 \norm{f}_{L^p(B(0, R))}^3.
  \end{multline*}
  We choose now a bounded domain $\Omega \subseteq \R^2 \times \R^2$
  such that the integral
  \begin{equation*}
    C_\Omega = \int_\Omega \frac{f(y) f(z)}
    {\abs{y} \abs{z} I_2\tonde{\abs{y} + \abs{z}, y + z}} \d y \d z
  \end{equation*}
  is finite and not zero
  (here $I_2$ is the function defined in lemma~\ref{lem:7}).
  This is possible when $f$ is not trivial,
  since $f(x)/\abs{x}$ is locally integrable and $1/I_2$ is bounded
  on compact subsets of $\mC_{+++}$.
  We divide both sides of the equation~\eqref{eq:57} by 
  the quantity $\abs{y} \abs{z} I_2\tonde{\abs{y} + \abs{z}, y + z}$
  and integrate with respect to $(y, z) \in \Omega$; we obtain
  \begin{equation} \label{eq:65}
    f(x) C_\Omega =
    \int_D F\tonde{\abs{x} + t, x + v} \varphi_\Omega(t, v) \d v \d t,
  \end{equation}
  for almost every $x \in \R^2$,
  where 
  \begin{equation*}
    \varphi_\Omega(t, v) = \1{I_2(t, v)}
      \int_\Omega \ddirac{t - \abs{y} - \abs{z} \\ v - y - z}
      \frac{\d y\d z}{\abs{y} \abs{z}} \le 1
  \end{equation*}
  is a bounded continuous function
  and the region $D$ is its support.
  The continuity of $f$ now follows from 
  the continuity of the right hand side in~\eqref{eq:65}
  by lemma~\ref{lem:15}.

  Suppose now that $f \in L^1_{\rm loc}(\R^2)$,
  then for $p > 2$ the functions $g = \abs{f}^{1/p} \in L^p_{\rm loc}(\R^2)$,
  $G = \abs{F}^{1/p}$ also solve equation~\eqref{eq:57}
  and it follows from the previous argument that $g$ is continuous.
  Hence, $\abs{f}$ is also continuous
  and so we have that $f \in L^p_{\rm loc}(\R^2)$ for any $p$.

  The continuity of $F$ comes easily from the equation 
  and the continuity of $f$:
  if~$t \ge r \ge 0$ and $\omega$ is a unit vector, we have
  \begin{equation*}
    F(t, r \omega) = f(0)
    f\TOnde{\frac{r + t}{2} \omega}
    f\TOnde{\frac{r - t}{2} \omega}.
  \end{equation*}
\end{proof}

\begin{lemma} \label{lem:16}
  If $f$ and $F$ are continuous functions which 
  solve equation~\eqref{eq:57}
  and~$f$ vanishes at one point
  then $f$ and $F$ vanish everywhere.
\end{lemma}

\begin{proof}
  Equation~\eqref{eq:57} implies that
  \begin{equation} \label{eq:66}
    f\tonde{\frac{x}3}^3 = F\tonde{\abs{x}, x} = f(x) f(0)^2.
  \end{equation}
  Suppose $f(x_0) = 0$ then $f(x_1) = 0$ for $x_1 = x_0 / 3$.
  By iterating this argument, $x_{k+1} = x_k/3$,
  we can construct a sequence of points $x_n$
  such that $f(x_n) = 0$ and $\lim_n x_n = 0$.
  By continuity it follows that $f(0) = 0$
  and by~\eqref{eq:66} $f$ must vanishes everywhere.
\end{proof}

\begin{lemma} \label{lem:17}
  Let $n \ge 1$.
  Let $\mN = \graffe{(t, x) \in \R \times \R^n: t = \abs{x}}$
  be the cone of future null vectors
  and $\mC = \graffe{(t, x) \in \R \times \R^n: t > \abs{x}}$
  the cone of future time-like vectors.
  Observe that $\mN + \mN = \ovl{\mC} = \mN \union \mC$.
  If $F: \mN \union \mC \to \C$ is a continuous solution
  of the conditional functional equation 
  \begin{equation} \label{eq:67}
    U, V \in \mN \implies F(U) F(V) = F(U + V)
  \end{equation}
  then $F$ is also a solution of the unconditional functional equation
  \begin{equation*}
    F(X) F(Y) = F(X + Y), \qquad \forall X, Y \in \mC.
  \end{equation*}
\end{lemma}

\begin{figure}
  \centering
  \begin{tikzpicture}[scale=0.2, >=latex]
    \draw[dotted, thin]
    (-12, 12) -- (1, -1) 
    (-1, -1) -- (12, 12);
    \draw[dashed, thin]
    (-11, 11) -- (-1, 20) -- (10, 10)
    (-8, 8) -- (-4, 12) -- (4, 4)
    (-3, 3) -- (3, 9) -- (6, 6)
    (-4, 12) -- (-1, 20) -- (3, 9);
    \draw[->, semithick]
    (0, 0) node[below] {$O$} -- (10, 10) node[below right] {$(a{+}c) U$};
    \draw[->, semithick] (0, 0) -- (6, 6) node[below right] {$a U$};
    \draw[->, semithick] (0, 0) -- (4, 4) node[below right] {$c U$};
    \draw[->, semithick] (0, 0) -- (-11, 11) node[below left] {$(b{+}d) V$};
    \draw[->, semithick] (0, 0) -- (-8, 8) node[below left] {$d V$};
    \draw[->, semithick] (0, 0) -- (-3, 3) node[below left] {$b V$};
    \draw[->, thick] (0, 0) -- (3, 9) node[above right] {$X$};
    \draw[->, thick] (0, 0) -- (-4, 12) node[above left] {$Y$};
    \draw[->, thick] (0, 0) -- (-1, 20) node[above] {$X{+}Y$};
  \end{tikzpicture}
  \caption{Construction for the proof of lemma~\ref{lem:17}.}
  \label{fig:3}
\end{figure}
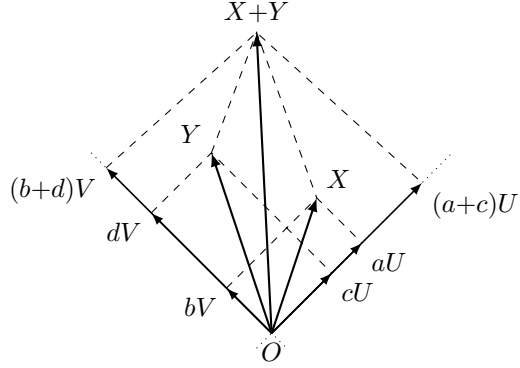

\begin{proof}
  Let $X$ and $Y$ be two vectors in $\mN \union \mC$
  which are not both in $\mN$.
  Let $\Pi$ be a two dimensional plane through the origin
  which contains $X$ and $Y$;
  the intersection of the plane $\Pi$ with the cone $\mN$ 
  is the union of two null directed half lines,
  \begin{equation*}
    \Pi \intersection \mN = (\R_+ U) \union (\R_+ V),
  \end{equation*}
  where $U$ and $V$ are two linearly independent vectors in $\mN$.
  We write $X$ and $Y$ as linear combinations of $U$ and $V$,
  \begin{equation*}
    X = a U + b V, \qquad
    Y = c U + d V,
  \end{equation*}
  for some non negative coefficients $a, b, c, d$.
  Then, using equation~\eqref{eq:67},
  \begin{multline*}
    F(X + Y) = F\Tonde{(aU + bV) + (cU +dV)} =
    F\Tonde{(a+c) U + (b+d) V} = \\
    = F\Tonde{(a+c) U} F\Tonde{(b+d) V} =
    \Tonde{F(aU) F(cU)} \Tonde{F(bV) F(dV)} = \\
    = \Tonde{F(aU) F(bV)} \Tonde{F(cU) F(dV)} = 
    F(aU + bV) F(cU + dV) = F(X) F(Y).
  \end{multline*}
\end{proof}

\begin{proposition} \label{pro:5}
  If $f: \R^2 \to \C$ and $F: \ovl{\mC_{+++}} \to \C$ 
  are non trivial locally integrable functions
  which satisfy the functional equation~\eqref{eq:57}
  then there exists constants $A \in \C$, $b \in \C^2$, $C \in \C$
  such that 
  \begin{equation*}
    f(x) = \exp\Tonde{A \abs{x} + b \cdot x + C}, \qquad
    F(t, x) = \exp\tonde{A t +  b \cdot x + 3C},
  \end{equation*}
  for (almost) all $(t, x) \in \mC_{+++}$.
\end{proposition}

\begin{proof}
  By lemma~\ref{lem:14},
  we may assume that $f$ and $F$ are continuous.
  By lemma~\ref{lem:16},
  we may assume that $f$ and $F$ never vanishes.
  Setting $y = 0$ and $z = 0$ in~\eqref{eq:57} we obtain
  $F\tonde{\abs{x}, x} = f(x) f(0)^2$.
  We define $G(t, x) = F(t, x) / F(0, 0)$;
  then
  \begin{equation*}
    G\Tonde{\abs{x}, x} = \frac{F\Tonde{\abs{x}, x}}{F(0, 0)} = 
    \frac{f(x)}{f(0)}, \qquad x \in \R^2.
  \end{equation*}
  We also have
  \begin{equation*}
    G\Tonde{\abs{x}, x} G\Tonde{\abs{y}, y} =
    \frac{f(x) f(y) f(0)}{f(0)^3} =
    \frac{F\Tonde{\abs{x} + \abs{y}, x + y}}{F(0, 0)} =
    G\Tonde{\abs{x} + \abs{y}, x + y}.
  \end{equation*}
  We can apply first lemma~\ref{lem:17} and then lemma~\ref{lem:10}
  to the function $G$ and obtain that $G(t, x) = \exp\tonde{A t + b \cdot x}$
  for some constants $A \in \C$ and $b \in C^2$.
  The result then follows by choosing $C$ so that $F(0, 0) = \exp(3C)$.
\end{proof}

\subsection{The equation~\eqref{eq:58}.}

\begin{lemma} \label{lem:18}
  It is possible to construct an open set $\Omega \subset \R^3 \times \R^3$,
  whose sections \mbox{$\Omega_x = \graffe{y: (x, y) \in \Omega}$}
  are dense in $\R^3$ for every $x \in \R^3$,
  and a pair of smooth maps $P, Q: \Omega \to \R^3$
  such that, for every $(x, y) \in \Omega$,
  \begin{gather}
    \label{eq:69}
    \abs{P(x, y)} + \abs{Q(x, y)} = \abs{x} + \abs{y}, \\
    \label{eq:70}
    P(x, y) + Q(x, y) = x + y, \\
    \notag
    \det\abs{\pder{P}{y}(x, y)} \ne 0, \qquad
    \det\abs{\pder{Q}{y}(x, y)} \ne 0.
  \end{gather}
\end{lemma}

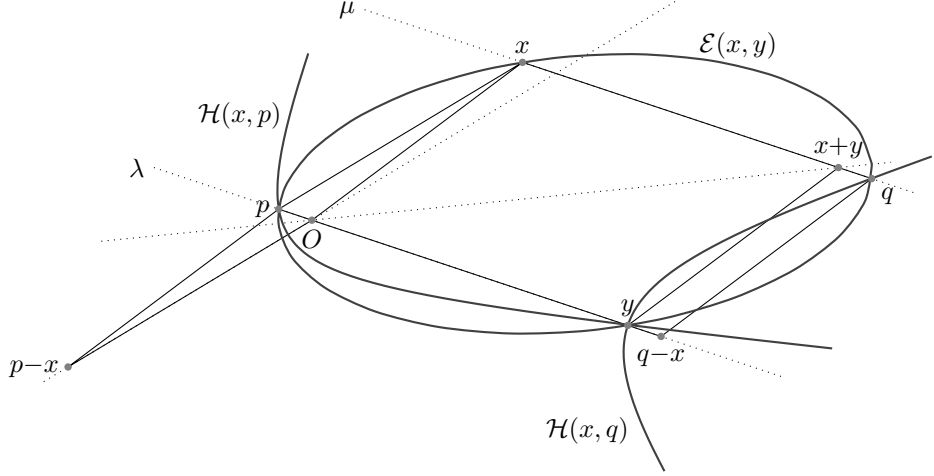
\begin{figure}
  \centering
  \begin{tikzpicture}[scale=0.7, >=latex]
    \draw[thin, dotted]
    (-4, -0.4) -- (11, 1.1)
    (-5.094, -3.069) -- (6.946, 4.185)
    (1, 4) node[left] {$\mu$} -- (11.5, 0.5)
    (-3, 1) node[left] {$\lambda$} -- (9, -3);
    \draw
    (6, -2) -- (10, 1)
    (4, 3) -- (0, 0) -- (-4.631, -2.790) -- (-0.631, 0.210) --
    (4, 3) -- (10.631, 0.790) -- (6.631, -2.210) -- (-0.631, 0.210);
    \draw[thick, black!75, parametric, smooth]
    plot[domain=0:6.283, id=Exy]
    function{5+5.634*cos(t)-0.259*sin(t), 0.5+0.563*cos(t)+2.597*sin(t)}
    plot[domain=-1.3:2.2, id=Hpx] function{%
      -2.315+1.857*cosh(t)+0.833*sinh(t), -1.385+1.119*cosh(t)-1.383*sinh(t)}
    plot[domain=-1:1.6, id=Hqx] function{%
      3.315+2.684*cosh(t)+0.648*sinh(t), -1.105-0.894*cosh(t)+1.945*sinh(t)};
    \draw
    (8.1, 2.9) node[above] {$\mE(x,y)$}
    (-0.4, 2) node[left] {$\mH(x,p)$}
    (6.2, -4) node[left] {$\mH(x,q)$};
    \fill[black!50, text=black]
    (0, 0) circle (2pt) node[below] {$O$}
    (4, 3) circle (2pt) node[above] {$x$}
    (6, -2) circle (2pt) node[above] {$y$}
    (10, 1) circle (2pt) node[above] {$x{+}y$}
    (-0.631, 0.210) circle (2pt) node[left] {$p$}
    (10.631, 0.790) circle (2pt) node[below right] {$q$}
    (-4.631, -2.790) circle (2pt) node[left] {$p{-}x$}
    (6.631, -2.210) circle (2pt) node[below] {$q{-}x$};    
  \end{tikzpicture}
  \caption{Constructions of the functions $P$, $Q$ and their inverses
    as described in lemma~\ref{lem:18}.}
  \label{fig:4}
\end{figure}

\begin{proof}
  The set $\Omega = \graffe{(x, y) \in \R^3 \times \R^3: x \times y \ne 0}$
  of linearly independent pair of vectors
  clearly has sections $\Omega_x$ dense in $\R^3$ for every $x$.
  Given $(x, y) \in \Omega$, the ellipsoid of revolution
  \begin{equation*}
    \mE(x, y) =
    \graffe{u \in \R^3: \abs{u} + \abs{x + y - u} = \abs{x} + \abs{y}},
  \end{equation*}
  with foci at $0$ and $x + y$ and which contains the points $x$ and $y$,
  is non degenerate and any line passing through one of the foci
  intersects the ellipsoid in exactly two points.
  In particular, the line $\lambda$ passing through $y$ and $0$
  intersects $\mE(x, y)$ in $y$ and in another point $p$;
  similarly, the line $\mu$ passing through $x$ and $x + y$
  intersects $\mE(x, y)$ in $x$ and in another point $q$.
  By symmetry we have $p + q = x + y$
  and from the definition of $\mE$ it follows that
  \begin{equation*}
    \abs{p} + \abs{q} = \abs{p} + \abs{x + y - p} = \abs{x} + \abs{y}.
  \end{equation*}
  It is evident from the geometric construction
  that the correspondence \mbox{$(x, y) \mapsto (p, q)$} is a smooth map
  as long as the vectors $x$ and $y$ remain linearly independent;
  moreover, when $x$ and $y$ are linearly independent
  we also have that
  $(x, p)$ and $(x, q)$ are pairs of linearly independent vectors.
  Setting $P(x, y) = p$ and $Q(x, y) = q$,
  we obtain two smooth maps $P, Q: \Omega \to \R^3$
  which satisfy~\eqref{eq:69} and~\eqref{eq:70}.

  To verify that, for fixed $x \in \R^3$,
  the maps $y \mapsto P(x, y)$ and $y \mapsto Q(x, y)$
  are locally invertible
  we provide a smooth geometric construction of their inverses.

  Given a pair of points $(x, p) \in \Omega$,
  we define $\mH(x, p)$ to be the branch of the hyperboloid
  with foci at $0$ and $p - x$ passing through the point $p$,
  \begin{equation*}
    \mH(x, p) =
    \graffe{u \in \R^3: \abs{u} - \abs{p - x - u} = \abs{p} - \abs{x}}
  \end{equation*}
  and we notice that it is non-degenerate
  since $p$ does not belong to the line passing through $0$ and $p - x$.
  The line passing through $p$ and $0$
  intersects $\mH(x, p)$ in $p$ and in another point $y_*$.
  The map $(x, p) \mapsto y_*$ 
  is smooth as long as $x$ and $p$ remain linearly independent.
  We claim that $P(x, y_*) = p$; indeed,
  $p$ belongs to the line passing through $y_*$ and $0$
  and from the definition of $\mH(x, p)$ it follows that
  \begin{equation*}
    \abs{p} + \abs{p - x - y_*} = \abs{x} + \abs{y_*},
  \end{equation*}
  which means that $p \in \mE(x, y_*)$.

  Similarly, given a pair of points $(x, q) \in \Omega$,
  we consider $\mH(x, q)$,
  the branch of the hyperboloid
  with foci at $0$ and $q - x$ passing through the point $q$.
  The line passing through $q - x$ and $0$
  intersects $\mH(x, q)$ in one point $y_{**}$,
  the vertex of the hyperboloid.
  The map $(x, q) \mapsto y_{**}$ 
  is smooth as long as $x$ and $q$ remain linearly independent
  and it easy to check that $Q(x, y_{**}) = q$.
  Indeed,
  since $q - x$ belongs to the line passing through $y_{**}$ and $0$,
  by a translation we have that
  $q$ belongs to the line passing through $x + y_{**}$ and $x$,
  moreover from the definition of $\mH(x, q)$ it follows that
  \begin{equation*}
    \abs{q} + \abs{q - x - y_{**}} = \abs{x} + \abs{y_{**}},
  \end{equation*}
  which means that $q \in \mE(x, y_{**})$.
\end{proof}

\begin{remark}
  Explicit formulae for the functions $P$ and $Q$
  constructed in the previous lemma are given by
  \begin{equation*}
    P(x, y) = \TOnde{\frac{x \cdot y - \abs{x} \abs{y}}
      {x \cdot y + \abs{x} \abs{y} + 2 \abs{y}^2}} y, \qquad
    Q(x, y) = x + y - P(x, y).
  \end{equation*}
\end{remark}

As in section~\ref{sec:case-w-n-3}
we let $\mC_{++} = \graffe{(t, v) \in \R \times \R^3: t > \abs{v}^2}$.

\begin{lemma} \label{lem:20}
  Let $f: \R^3 \to \C$ and $F: \ovl{\mC_{++}} \to \C$
  be functions which solve equation~\eqref{eq:58}.
  If $f$ is locally integrable then $f$ and $F$ are continuous functions.
\end{lemma}

\begin{proof}
  Let $P$ and $Q$ be the functions constructed in lemma~\ref{lem:18}.
  If $f$ and $F$ are solutions to~\eqref{eq:58}
  it follows that
  \begin{equation*}
    f(x) f(y) = f\Tonde{P(x, y)} f\Tonde{Q(x, y)}, \quad
   \text{for a.e. $(x, y) \in \R^3 \times \R^3$},
  \end{equation*}
  and the lemma then becomes a corollary of proposition~\ref{pro:2}.
\end{proof}

Once the continuity of locally integrable solutions to~\eqref{eq:58}
is established, one then proceeds in the same manner as in the
previous subsection and obtains the following result.

\begin{proposition} \label{pro:4}
  If $f: \R^3 \to \C$ and $F: \ovl{\mC_{++}} \to \C$ 
  are non trivial locally integrable functions
  which satisfy the functional equation
  \begin{equation*}
    f(x) f(y) = F\Tonde{\abs{x} + \abs{y}, x + y},
  \end{equation*}
  for all $x, y \in \R^3$,
  then there exists constants $A \in \C$, $b \in \C^3$, $C \in \C$
  such that 
  \begin{equation*}
    f(x) = \exp\Tonde{A \abs{x} + b \cdot x + C}, \qquad
    F(t, x) = \exp\tonde{A t +  b \cdot x + 2 C},
  \end{equation*}
  for (almost) all $(t, x) \in \mC_{++}$.
\end{proposition}

\section*{%
  Acknowledgement}

I am grateful to S.\ Klainerman
for interesting discussions and for suggesting the proof of lemma~\ref{lem:17};
to A.\ J\'arai 
for comments on regularity properties of solutions to functional equations;
and to the anonymous referee
for suggesting a simpler way of computing
the integrals of lemma~\ref{lem:3} and lemma~\ref{lem:7}
based on Lorentz invariance,
for pointing out a flaw in the first version of the proof of lemma~\ref{lem:8},
and for suggesting the example described in remark~\ref{rem:6}.

\bibliographystyle{amsplain}

\end{document}